\documentclass[11pt]{article} 
\usepackage{latexsym} 
\input{amssym.def} 
\input{amssym} 
\newfont{\fra}{eufm10 scaled 1095} 
\newfont{\Bb}{msbm10 scaled 1095} 
\newfont{\Bbg}{msbm10 scaled 1280} 
 
\newcommand\RR{{\mbox{\Bb R}}}

\newcommand\X{\mbox{\fra X}} 
\newcommand\fg{{\frak{g}}} 
\newcommand\fh{{\frak h}} 
\newcommand\fri{{\frak i}} 
 
\newcommand\fl{{\frak l}}

\newcommand\fq{{\frak q}} 
\newcommand\fa{{\frak a}} 
 
\newcommand\fd{{\frak d}}

\newcommand\fz{{\frak z}}

\newcommand\cZ{{\cal Z}} 
\newcommand\cH{{\cal H}}

\newcommand\ph{\varphi}

\newcommand{\gl}{\mathop{{\frak g \frak l}}}

\newcommand{\Aut}{\mathop{{\rm Aut}}} 
\newcommand{\GL}{\mathop{{\it GL}}}

\newcommand{\Id}{\mathop{{\rm id}}} 
 
\newcommand{\ad}{\mathop{{\rm ad}}}

\newcommand{\sgn}{\mathop{{\rm sgn}}}

\newcommand{\diag}{\mathop{{\rm diag}}} 
 
\newcommand{\Span}{\mathop{{\rm span}}}

\newcommand{\cZQ}{{\cZ^{2}_{Q}(\fl,\fa)}} 
\newcommand{\cHQ}{{\cH^{2}_{Q}(\fl,\fa)}} 
\newcommand\ip{\mbox{$\langle\cdot \,,\cdot \rangle$}} 
\newcommand\ipa{{\langle\cdot \,,\cdot \rangle_\fa}} 
\newcommand\lb{{[\cdot\,,\cdot]}} 
 
\newcommand\dd{\fd_{\alpha,\gamma}(\fl,\fa)} 
\newcommand\proof{{\sl Proof. }} 
\newcommand{\qed}{\hspace*{\fill}\hbox{$\Box$}\vspace{2ex}} 
\newtheorem{theo}{Theorem} 
\newtheorem{pr}{Proposition} 
\newtheorem{de}{Definition} 

\newtheorem{lm}{Lemma} 
\begin{document} 
\title{Nilpotent metric Lie algebras of small dimension} 
\author{Ines Kath} 
\maketitle 
\begin{abstract} \noindent In \cite{KO2} we developed a general classification 
scheme for metric Lie algebras, i.e.~for finite-dimensional Lie algebras 
equipped with a non-degenerate invariant inner product. Here we determine all 
nilpotent Lie algebras $\fl$ with $\dim \fl'=2$ which are used in this scheme. 
Furthermore, we classify all nilpotent metric Lie algebras of dimension at most 
10.  
\end{abstract} 
\tableofcontents 

\section{Introduction} 

In \cite{KO2} we developed a structure theory for metric Lie algebras, 
i.e.~for Lie algebras with invariant non-degenerate inner product or, 
equivalently, for simply-connected Lie groups with a bi-invariant 
pseudo-Riemannian 
metric. We used 
this structure theory in order to give a description of the moduli space of all 
isomorphism 
classes of indecomposable non-simple metric Lie algebras as 
\begin{equation}\label{x} 
\coprod_{(\frak l,\frak a)}{\cal H}^{2}_{Q}(\frak l,\frak a)_{0}/G_{(\frak 
l,\frak 
a)}, 
\end{equation} 
where the union is taken over all isomorphism classes of pairs $(\frak l,\frak 
a)$ of Lie algebras $\frak l$ and semi-simple orthogonal $\frak l$-modules 
$\frak a$. 
Here ${\cal H}^{2}_{Q}(\frak l,\frak a)_{0}$ denotes a certain subset of 
the second quadratic cohomology ${\cal 
H}^{2}_{Q}(\frak l,\frak a)$ (see also Sections \ref{S2} and \ref{S4} for a 
definition of 
these sets) and $G_{(\frak l,\frak a)}$ is the automorphism group of 
the pair $(\frak l,\frak a)$. 

Moreover, in \cite{KO2} we gave explicitly the map which assigns an element of  
(\ref{x}) to 
each isomorphism class of indecomposable metric Lie algebras as well as its 
inverse 
map. The 
construction of these maps relies on the fact that for each metric Lie algebra 
without 
simple ideals there is a canonical isotropic ideal $\fri(\fg)\subset \fg$ such 
that 
$\fa:= 
\fri(\fg)^\perp/\fri(\fg)$ is abelian. 

The description (\ref{x}) of the moduli 
space of isomorphism classes of metric Lie algebras allows a 
systematic approach to the construction and classification of metric 
Lie algebras. Of course it is far from being an 
explicit classification (e.g.~a list). A full classification would require that 
we 
can determine all Lie algebras $\fl$ for which 
${\cal H}^{2}_{Q}(\frak l,\frak a)_{0}$ is not empty for some 
orthogonal $\fl$-module $\fa$. These Lie algebras are called 
admissible. However, although admissibility is a strong condition 
it seems to be hard to give a classification of these Lie algebras. 
Another problem is the explicit computation of the cohomology sets 
which includes for example the classification of $GL(l,\RR)$-orbits of 3-forms 
on 
an $l$-dimensional vector space. Such a 
classification is known only for $l\le9$. 

However, (\ref{x}) yields a general classification scheme which can be 
used to obtain a full classification 
for metric Lie algebras satisfying suitable additional assumptions. Such 
assumptions can be, e.g., restrictions on the index of the inner product or on 
the structure of the Lie algebra. These restrictions give additional conditions 
for the Lie algebras $\fl$ occuring in (\ref{x}). Hence, in order to get a 
classification from (\ref{x}) one has 
first to determine all admissible Lie algebras $\fl$ which satisfy these 
additional conditions and afterwards one has to determine orbit sets of 
cohomology classes of these Lie algebras. For example, the classification of 
metric Lie algebras with index $p$ leads to the classification problem for 
admissible Lie algebras of dimension $\dim \fl\le p$. In \cite{KO1} and 
\cite{KO2} we show how one can solve this problem for small $p$. In particular, 
we give a classification of all metric Lie algebras whose invariant inner 
product is of index two or three.  

We see that the classification of admissible Lie algebras within a certain class 
is a main step in the solution of the 
original classification problem for metric Lie algebras (with additional 
properties). 
In general, the classification 
of admissible Lie algebras even within a 
certain class of Lie algebras seems to be complicated. However, often it is much 
easier than the determination  of all Lie algebras of this class.

Let us consider another suitable condition which allows to make (\ref{x}) more 
explicit. Namely, let us consider only indecomposable metric Lie algebras whose 
canonical isotropic ideal $\fri(\fg)$ is ``almost central'' (this 
means that the codimension of $\fz(\fg)\subset\fri(\fg)$ in $\fri(\fg)$ is 
small). The case $\fri(\fg)=\fz(\fg)$ has been studied in \cite{KO1}. In 
particular,  the general classification scheme has been specialised to the case 
of metric Lie algebras with maximal isotropic centre. If $\fri(\fg)=\fz(\fg)$, 
then we have to consider only abelian Lie algebras $\fl$ in (\ref{x}). All 
abelian Lie algebras are admissible. 
In the more general case  $\fz(\fg)\subset\fri(\fg)$ with small codimension we 
are led to the investigation of admissible Lie algebras $\fl$ with small 
nilpotent radical $R(\fl)$. By definition $R(\fl)$ is the minimal ideal such 
that 
the adjoint representation of $\fl$ on 
$\fl/R(\fl)$ is semi-simple, e.g.~$R(\fl)=\fl'$ for nilpotent~$\fl$. 

In the first part of this paper we solve the classification problem for          
        nilpotent admissible Lie algebras $\fl$
whose nilpotent radical $R(\fl)=\fl'$ is two-dimensional. We can prove that such 
Lie algebras are direct sums $\fg\oplus\RR^k$, where $\fg$ is nilpotent and 
admissible of dimension at most 6. The precise classification result 
is stated in Section~\ref{S3}, Proposition~\ref{P2}. Solvable non-nilpotent 
admissible Lie algebras with two-dimensional nilpotent radical and 
solvable admissible Lie algebras with one-dimensional nilpotent radical 
were already classified in \cite{KO2}. 

In the second part we apply the general classification 
scheme for metric Lie algebras to low-dimensional nilpotent metric 
Lie algebras. We use the classification results from the first 
part of the paper to determine all nilpotent metric Lie 
algebras of dimension~${\le 10}$, see Theorem \ref{T1} at the end of 
this paper.

\section{Admissible cohomology classes} 
\label{S2}
In \cite{KO2} we defined the quadratic cohomology $\cHQ$ for a Lie 
algebra $\fl$ and an orthogonal $\fl$-module $\fa$. Let us recall 
this definition.
 An {\it orthogonal $\fl$-module} is a 
tuple $(\rho,\fa,\ip_{\fa})$ (also $\fa$ or $(\rho,\fa)$ in abbreviated 
notation) 
consisting of a finite-dimensional pseudo-Euclidean vector space 
and a representation $\rho$ of $\fl$ on $\fa$ satisfying 
$$\langle \rho(L)A_1,A_2\rangle_\fa + \langle A_1, \rho(L)A_2\rangle_\fa =0$$
for all $L\in\fl$ and $A_1,A_2\in \fa$. 

For $\fl$ and (any $\fl$-module) $\fa$ we have the
standard cochain complex $(C^*(\fl,\fa),d)$ and corresponding cohomology groups
$H^p(\fl,\fa)$. 
If $\fa$ is the one-dimensional trivial
representation, then we denote this cochain complex also by  $C^*(\fl)$. 

We define the product
$$\langle \cdot \wedge \cdot \rangle :\ C^p(\fl,\fa)\times 
C^q(\fl,\fa)\longrightarrow C^{p+q}(\fl)$$
by the composition
$$C^p(\fl,\fa)\times C^q(\fl,\fa)\stackrel{\wedge}{\longrightarrow}
C^{p+q}(\fl,\fa\otimes \fa) \stackrel{\ip_\fa}{\longrightarrow} C^p(\fl).$$

Let $p$ be even. Then the group of quadratic $(p-1)$-cochains is the group
$${\cal C} ^{p-1}_Q(\fl,\fa)=C^{p-1}(\fl,\fa)\oplus C^{2p-2}(\fl) $$
with  group operation defined by
$$ (\tau_1,\sigma_1)*(\tau_2,\sigma_2)=(\tau_1+\tau_2, \sigma_1
+\sigma_2 +\textstyle{\frac12} \langle \tau_1\wedge \tau_2\rangle)\,.$$
Now we consider the set  
$${\cal Z} ^{p}_Q(\fl,\fa)=\{(\alpha,\gamma) \in C^{p}(\fl,\fa)\oplus
C^{2p-1}(\fl) \mid d\alpha=0,\
d\gamma=\textstyle{\frac12}\langle\alpha \wedge\alpha\rangle\} $$
of so-called quadratic $p$-cocycles. The group  
${\cal C} ^{p-1}_Q(\fl,\fa)$ acts on ${\cal Z} ^{p}_Q(\fl,\fa)$ by
$$(\alpha,\gamma)(\tau,\sigma)=\Big(\,\alpha +d\tau,\gamma +d\sigma
+\langle(\alpha +\textstyle{\frac12} d\tau)\wedge\tau\rangle\,\Big).$$
and we define the quadratic cohomology set ${\cal H}
^{p}_Q(\fl,\fa)~:={\cal Z} ^{p}_Q(\fl,\fa)/ {\cal C}
^{p-1}_Q(\fl,\fa)$. As usual, we denote the equivalence class of
$(\alpha,\gamma)\in {\cal Z} ^{p}_Q(\fl,\fa)$ in ${\cal H}
^{p}_Q(\fl,\fa)$ by $[\alpha,\gamma]$.

Let us now recall the 
definition of admissible cohomology 
classes from \cite{KO2}. In general, admissible cohomology classes are certain 
elements of $\cHQ$ for a Lie algebra $\fl$ and a semi-simple orthogonal 
$\fl$-mod\-ule $\fa$. Here we will give the definition of admissibility only for 
nilpotent Lie algebras. So we have the two following simplifications compared to 
\cite{KO2}: 

If $\fl$ is 
a nilpotent Lie algebra, then 
$H^{*}(\fl,\fa)=H^{*}(\fl,\fa^{\fl})$ holds for any semi-simple $\fl$-mod\-ule 
$\fa$ (see \cite{D}). This implies 
that also $\cHQ=\cH_Q^{2}(\fl,\fa^{\fl})$ holds for any 
orthogonal semi-simple $\fl$-module $\fa$. 

For a Lie algebra $\fl$ we denote by $\fl^1=\fl,\dots, 
\fl^k=[\fl,\fl^{k-1}],\dots$ 
the lower central series. As usual we often denote $\fl^2$ also by $\fl'$. If 
$\fl$ 
is nilpotent, then its $k$-th nilpotent radical  $R_k(\fl)$ equals $\fl^{k+1}$. 
\begin{de} 
Let $\fl$ be a nilpotent Lie algebra and let $(\rho,\fa,\ipa)$ be a 
semi-simple 
orthogonal $\fl$-module. Let 
$m$ be such 
that $\fl^{m+2}=0$. Put $\fl_{(0)}=\fz(\fl)\cap \ker \rho$ and 
$\fl_{(k)}=\fz(\fl)\cap 
\fl^{k+1}$ for $k\ge1$. Take a cohomology class in $\cHQ$ and represent it by a 
cocycle $(\alpha,\gamma)$ satisfying $\alpha(\fl,\fl)\subset\fa^\fl$. Then  
$[\alpha,\gamma]\in\cHQ$ is called 
admissible if and only if 
the following 
conditions $(A_k)$ and $(B_k)$ hold for all $0\le k\le m$. 
\begin{enumerate} 
\item[$(A_k)$] 
Let $L_0\in \fl_{(k)}$ be such that there exist 
elements 
$A_0\in \fa$ and $Z_0\in (\fl^{k+1})^*$ satisfying 
\begin{enumerate} 
\item[(i)] $\alpha(L,L_0)=0 $, 
\item[(ii)] $\gamma(L,L_0,\cdot)=-\langle A_0,\alpha(L,\cdot)\rangle_\fa 
+\langle 
Z_0, [L,\cdot]_\fl\rangle$ as an element of $(\fl^{k+1})^*$, 
\end{enumerate} 
for all $L\in\fl$, then $L_0=0$. 
\item[$(B_k)$] The subspace $\alpha(\ker \lb_{\fl\otimes \fl^{k+1}})\subset 
\fa$ is non-degenerate, where $\ker \lb_{\fl\otimes \fl^{k+1}}$ is the kernel of 
the 
map 
$\lb:\fl\otimes \fl^{k+1}\rightarrow \fl$. 
\end{enumerate} 
We denote the set of all admissible cohomology classes in $\cHQ$ by 
$\cHQ_{\sharp}$. 
A Lie algebra $\fl$ is called admissible if there is a semi-simple 
orthogonal $\fl$-module $\fa$ such that $\cHQ_{\sharp}\not=\emptyset$. 
\end{de} 

\section{Nilpotent admissible Lie algebras with 2-dimensional radical} 
\label{S3}
In the following we will often describe a Lie algebra by giving a 
basis and some of the Lie brackets. In this case we always assume that all other 
brackets of basis vectors vanish. If we do not mention the basis 
explicitly, then we assume that all basis vectors appear in one of 
the bracket relations (on the left or the right hand side). 

Using this convention we define 
\begin{eqnarray*} 
\fh(1) &=&\{ [X_1,X_2]=Y\}\\ 
\fg_{4,1} &=&\{ [X_1,Z]=Y,\ [X_1,X_2]=Z\},\\ 
\fg_{5,2} &=&\{ [X_1,X_2]=Y,\ [X_1,X_3]=Z\},\\ 
\fg_{6,4} &=&\{ [X_1,X_2]=Y,\ [X_1,X_3]=Z,\ [X_3,X_4]=Y\},\\ 
\fg_{6,5} &=&\{ [X_1,X_2]=Y,\ [X_1,X_3]=Z,\ [X_2,X_4]=Z,\ [X_3,X_4]=- 
Y\}. 
\end{eqnarray*} 
Note that $\fg_{6,4}$ and $\fg_{6,5}$ are not isomorphic since $\dim 
[X_2,\fg_{6,4}]=1$ but $\dim [L,\fg_{6,5}]=2$ for all $L\in\fg_{6,5}$. 

\begin{pr} \label{P1}
The Lie algebras $\fh(1)$, $\fg_{4,1}$, $\fg_{5,2}$, $\fg_{6,4}$, 
$\fg_{6,5}$ are 
admissible. 
\end{pr} 
\proof In \cite{KO2} we proved that $\fh(1)$ is admissible. In Propositions 
\ref{P5} 
and \ref{g41}  we will see that $\fg_{4,1}$ and $\fg_{5,2}$ are also admissible. 
Let 
us verify now that the statement holds for $\fg_{6,4}$ and 
$\fg_{6,5}$. First we consider $\fl=\fg_{6,4}$. Take $\fa=\RR^{2,2}$ and 
$\rho=0$. Let  
$A_1,A_2,A_3,A_4$ be a Witt basis of $\fa$, i.e. $\langle A_1,A_3\rangle 
=\langle A_2,A_4\rangle =1$ and $\langle A_i,A_j\rangle =0$ for the remaining 
pairs $1\le i\le j\le4$.  
We define $\alpha\in C^2(\fl,\fa)$ by 
\begin{eqnarray*} 
\alpha(X_1,Y)=-\alpha(X_4,Z)=A_1,&& \alpha(X_3,Y)=\alpha(X_2,Z)=A_2\\ 
\alpha(X_3,Z)=A_3,&& \alpha(X_1,Z)=A_4\\ 
\alpha(X_2,Y)=\alpha(X_4,Y)=0,&& \alpha(X_i,X_j)=\alpha(Y,Z)=0. 
\end{eqnarray*} 
Then it is easy to see that $d\alpha=0$ and $\langle \alpha\wedge\alpha\rangle 
=0$. 
Hence $(\alpha,0)\in\cZ_Q^2(\fl,\fa)$. Let us now show that $[\alpha,0]$ is 
admissible. Take $L_0=cY+dZ\in\fl_{(k)}\subset\Span\{Y,Z\}$. Then 
$\alpha(L_0,\fl)=0$ 
would imply $\alpha(cY+dZ,X_1)=-cA_1-dA_4=0$, hence $c=d=0$, and therefore 
$L_0=0$. 
Thus Condition $(A_k)$ is satisfied. It remains to check Conditions $(B_k)$ for 
$k=0,1$ since $\fl^3=0$. These are satisfied because of $\alpha(\ker 
\lb_{\fl\otimes 
\fl^{k+1}})=\fa$ for $k=0,1$. 

Now we consider $\fl=\fg_{6,5}$. Let $(\rho,\fa)$ be as above and 
define $\alpha\in C^2(\fl,\fa)$ by 
\begin{eqnarray*} 
&\alpha(X_1,Y)=\alpha(X_4,Z)=A_1,\qquad \alpha(X_3,Y)=\alpha(X_2,Z)=A_2&\\ 
&\alpha(X_2,Y)=-\alpha(X_3,Z)=A_3,\qquad \alpha(X_4,Y)=-\alpha(X_1,Z)=A_4&\\ 
& \alpha(X_i,X_j)=\alpha(Y,Z)=0. &
\end{eqnarray*} 
In the same way as above one verifies $[\alpha,0]\in \cHQ_{\sharp}$.
\qed 
\begin{pr}\label{P2} 
If $\fl$ is an admissible nilpotent Lie algebra with $\dim \fl'=2$, then 
$\fl$ is 
isomorphic to one of the (admissible) Lie algebras 
$$\fh(1) \oplus \fh(1) \oplus \RR^k,\ \fg_{4,1}\oplus \RR^k,\ 
\fg_{5,2}\oplus \RR^k,\ 
\fg_{6,4}\oplus \RR^k,\ \fg_{6,5}\oplus \RR^k.$$ 
\end{pr} 
\proof Let us verify that all these Lie algebras are admissible. First notice 
that $\RR^k$ is admissible. Indeed, let $X_1,\dots,X_k$ be a basis of $\RR^k$ 
and take $\fa$ and $\alpha\in C^2(\fl,\fa)$ such that $\alpha(X_i,X_j)=A_{ij}$ 
for $1\le i< j\le k$ and $\{A_{ij}\}_{1\le i< j\le k}$ is an orthonormal basis 
of $\fa$. Then $[\alpha,0]\in\cHQ_\sharp$. If we use now Proposition \ref{P1} 
and the fact that direct sums of admissible Lie algebras are admissible the 
assertion follows.

Now we prove that each admissible nilpotent Lie algebra $\fl$ with $\dim \fl'=2$ 
is 
isomorphic to one of the mentioned Lie algebras. We distinguish between two 
cases: $\fl'\not\subset \fz(\fl)$ (case~I) and 
$\fl'\subset \fz(\fl)$ (case~II). 

\subsubsection*{Case I: $\fl'\not\subset \fz(\fl)$} 
The representation of $\fl$ on $\fl'$ is nilpotent and non-trivial. Hence we 
may 
choose a basis $Y,Z$ of $\fl'$ and a vector $X_1$ in $\fl\setminus \fl'$ such 
that 
\begin{equation}\label{E1} 
[X_1,Z]=Y\,. 
\end{equation} 
In particular, since $\fl$ is nilpotent this implies 
\begin{equation}\label{E3} 
[\fl,Y]=0\,. 
\end{equation} 
Using this we can see that $[X_1,\fl]$ is not contained in $\RR\cdot Y$. 
Indeed, 
$[X_1,\fl]\subset\RR\cdot Y$ would imply 
$[X_1,\fl']=[X_1,[\fl,\fl]]=[[X_1,\fl],\fl]\subset [Y,\fl]=0$, which contradicts 
$Y=[X_1,Z]\in 
[X_1,\fl']$. We conclude that there is a vector $X_2\in\fl\setminus \fl'$ such 
that 
\begin{equation}\label{E2} 
[X_1,X_2]=Z\,. 
\end{equation} 
Here we may assume 
\begin{equation} 
[X_2,Z]=0\,.\label{E8} 
\end{equation} 
By (\ref{E1}) and (\ref{E2}) it is possible to choose a vector space 
decomposition 
$$\fl=\Span\{X_1,X_2\}\oplus V \oplus \fl' $$ 
of $\fl$ such that $[X_1,V]=0$ and $[X_2,V]\subset \RR\cdot Y$. In 
particular, this 
implies $[X_1, [V,V]]=[[X_1,V],V] =0$, hence 
$$[V,V]\subset \RR\cdot Y$$ 
by (\ref{E1}). Moreover, $[X_1,V]=0$ and $[X_2,V]\subset \RR\cdot Y$ together 
with 
(\ref{E2}) gives 
$$[V,Z]=[V,[X_1,X_2]]= [X_1,[V,X_2]]\subset [X_1,\RR\cdot Y]=0.$$ 
Now we distinguish between the 
cases $[X_2,V]\not=0$ and $[X_2,V]=0$. 

\subsubsection*{Case I.1: $[X_2,V]\not=0$} 
{\bf Claim.} A Lie algebra $\fl$ which satisfies the conditions of case I.1 
is not 
admissible.\\[2ex] 
\proof By (\ref{E1}) -- (\ref{E8}) and our choice of $V$ we find a basis 
$X_3,\dots,X_l$ 
of $V$ such that 
$$ 
\fl=\{\,[X_1,X_2]=Z,\,[X_1,Z]=Y,\,[X_2,X_3]=Y,\, 
[X_i,X_j]=y_{ij}Y,\,i,j\ge 3 
\,\}
$$ 
for suitable $y_{ij}\in\RR$. Assume that $\fl$ 
is admissible. Then we can choose a semi-simple orthogonal $\fl$-module $\fa$ 
and 
$[\alpha,\gamma]\in \cHQ$ such that 
$[\alpha,\gamma]$ is admissible. 
As explained above we may assume $\alpha(\fl,\fl)\subset \fa^\fl$. Hence 
$d\alpha=0$ 
implies 
\begin{eqnarray} 
0&=&\alpha([X_2,X_3],Z)\,=\,\alpha(Y,Z) \label{E4}\\ 
0&=&\alpha([X_1,X_2],X_3)+\alpha([X_2,X_3],X_1)\,=\,\alpha(Z,X_3)+ 
\alpha(Y,X_1)\label{E5} 
\\ 
0&=&\alpha([X_1,Z],X_j)\,=\,\alpha(Y,X_j),\quad 
j\ge2\,.\label{E6} 
\end{eqnarray} 
Because of $ \langle \alpha\wedge \alpha \rangle=2 d\gamma$ we have 
\begin{eqnarray*} 
\lefteqn{\langle\alpha(X_1,X_3),\alpha(Y,Z)\rangle + 
\langle\alpha(X_3,Y),\alpha(X_1,Z)\rangle + 
\langle\alpha(Y,X_1),\alpha(X_3,Z)\rangle\qquad}\\ 
&& \qquad=d\gamma(X_1,X_3,Y,Z)=-\gamma([X_1,Z],X_3,Y)= -\gamma(Y,X_3,Y)=0 
\end{eqnarray*} 
and by (\ref{E4}) -- (\ref{E6}) this yields 
$\langle\alpha(Y,X_1),\alpha(Y,X_1)\rangle=0\,.$ 
Summarizing we obtain 
\begin{equation} \label{E7} 
\alpha(Y,Z)=0,\quad \langle\alpha(Y,X_1),\alpha(Y,X_1)\rangle=0,\quad 
\alpha(Y,X_j)=0,\ j\ge2\,. 
\end{equation} 
Now let us consider Condition $(B_2)$. Since $\fl^3=\RR\cdot Y\subset\fz(\fl)$ 
it is 
satisfied if 
and only if the space $\alpha(\fl,Y)$ is non-degenerate. Now (\ref{E7}) 
implies that 
$(B_2)$ is satisfied if and only if $\alpha(Y,\fl)=0$. But if $\alpha(Y,\fl)=0$, 
then Condition $(A_2)$ is not satisfied. Indeed, $L_0=Y\not=0$, $A_0=0$, $Z_0=0$ 
obviously 
satisfy $(A_2)\,(i)$ and $(A_2)\,(ii)$ since $\fl^3=\RR\cdot Y$ 
is one-dimensional. Thus we obtain a contradiction and $\fl$ is not admissible. 
\qed 

\subsubsection*{Case I.2: $[X_2,V]=0$} 

{\bf Claim.} An admissible Lie algebra $\fl$ which satisfies the conditions 
of case I.2 is isomorphic to $\fg_{4,1}\oplus \RR^k$.\\[2ex] 
\proof 
By (\ref{E1}) -- (\ref{E8}) and our choice of $V$ we find a basis 
$X_3,\dots,X_l$ of $V$ 
such that 
$$ 
\fl=\{\, 
[X_1,X_2]=Z,\,[X_1,Z]=Y,\, 
[X_i,X_j]=y_{ij}Y,\,i,j\ge 3 
\,\}\,. 
$$ 
for suitable $y_{ij}\in\RR$, $i,j\ge 3$. 
Suppose $[\alpha,\gamma]\in \cHQ$ is admissible and $\alpha(\fl,\fl)\subset 
\fa^\fl$. The cocycle conditions 
$$d\alpha(X_1,X_2,Y)=0,\ d\alpha(X_i,X_j,X_1)=0,\ d\alpha(X_1,X_2,Z)=0,\ 
d\alpha(X_1,X_i,Z)=0$$ 
for $i,j\ge 3$ yield 
$$ \alpha(Y,Z)=0,\ y_{ij}\alpha(Y,X_1)=0, \ \alpha(Y,X_2)=0, \ 
\alpha(Y,X_i)=0,\,i,j\ge 3,$$ 
respectively. Assume that $y_{ij}\not=0$ for some $i,j\ge 3$. Then 
$\alpha(Y,\fl)=0$ follows. 
In this case $L_0= Y\in \fz(\fl)\cap\fl^3$, $A_0=0$ satisfy 
$(A_2)\,(i)$ and 
$(A_2)\,(ii)$ since $\fl^3=\RR\cdot Y$ is one-dimensional. Since $Y\not=0$ we 
see that $(A_2)$ is not satisfied, a contradiction. Consequently, $y_{ij}=0$ for 
all $i,j\ge 3$, which proves the claim. 
\qed 

\subsubsection*{Case II: $\fl'\subset \fz(\fl)$} 
\begin{lm}\label{L1} 
There exists a 3-dimensional subspace $\bar \fl$ of $\fl$ such that 
$[\,\bar\fl,\bar\fl\,]=\fl'$. Moreover, we can choose a basis $X_1,X_2,X_3$ of 
$\bar 
\fl$ and a basis $Y,Z$ of $\fl'$ such that 
\begin{equation} \label{E10n}
[X_1,X_2]=Y,\ [X_1,X_3]=Z,\ [X_2,X_3]=0. 
\end{equation} 
\end{lm} 
\proof 
Since $\dim \fl'=2$ we can choose vectors $L_1,\dots,L_4$ such 
that 
$\fl'=\Span\{[L_1,L_2],$ $[L_3,L_4]\}$. We consider 
$\fl_1=\Span \{L_1,L_2\}$ and $\fl_2=\Span \{L_3,L_4\}$. If 
$[\fl_1,\fl_2]\not=0$, then 
we may assume $[L_1,L_3]\not=0$. In this case at least one of the pairs 
$[L_1,L_3], 
[L_1,L_2]$ and $[L_1,L_3], [L_3,L_4]$ consists of two linearly independent 
vectors and we 
can choose $\bar \fl$ correspondingly. If $[\fl_1,\fl_2]=0$, then e.g. 
$\bar 
\fl=\Span\{L_1, L_3,L_2+L_4\}$ satisfies $[\,\bar\fl,\bar\fl\,]=\fl'$. 

We can choose linearly independent vectors 
$X_1,X_2,X_3,Y,Z$ of $\fl$ 
such that $[X_1,X_2]=Y$ and $[X_1,X_3]=Z$. If $[X_2,X_3]=yY+zZ$, then $\bar 
X_2:=X_2-zX_1$ and $\bar X_3:=X_3+yX_1$ satisfy $[X_1,\bar X_2]=Y$, 
$[X_1,\bar X_3]=Z$ 
and $[\bar X_2, \bar X_3]=0$. 
\qed 

\begin{lm}\label{L3} Let $X_1,X_2,X_3,Y,Z$ be as in Lemma \ref{L1}. If 
$[\alpha,\gamma]\in\cH^2_Q(\fl,\fa)$ and $\alpha(\fl,\fl)\subset \fa^\fl$, then 
we 
have 
\begin{itemize} 
\item[(i)] $\alpha(Y,Z)=0$; 
\item[(ii)] $\alpha(Y,L)=0$ for all $L\in\fl$ satisfying $[L,X_1]=[L,X_2]=0$; 
\item[(iii)] $\alpha(Z,L)=0$ for all $L\in\fl$ satisfying $[L,X_1]=[L,X_3]=0$; 
\item[(iv)] $\langle\alpha(U_1,L_1),\alpha(U_2,L_2)\rangle = 
\langle\alpha(U_1,L_2),\alpha(U_2,L_1)\rangle$ for all $U_1,U_2\in \fl'$ and 
$L_1,L_2\in\fl$. 
\end{itemize} 
\end{lm} 
\proof 
Assertions {\it (i)}, {\it (ii)}, {\it (iii)} follow from the cocycle 
condition for 
$\alpha$, from $\fl'\subset 
\fz(\fl)$ and from the special conditions on $L\in\fl$ in {\it (ii)} and {\it 
(iii)}, 
respectively: 
$$\begin{array}{lclclcl} 
\alpha(Y,Z)&=&\alpha([X_1,X_2],Z)&=&\alpha([Z,X_2],X_1) + 
\alpha([X_1,Z],X_2) &=&0\\[1ex] 
\alpha(Y,L)&=&\alpha([X_1,X_2],L)&=&\alpha([L,X_2],X_1) + 
\alpha([X_1,L],X_2) 
&=&0\\[1ex] 
\alpha(Z,L)&=&\alpha([X_1,X_3],L)&=&\alpha([L,X_3],X_1) + 
\alpha([X_1,L],X_3) 
&=&0. 
\end{array}$$

As for assertion {\it (iv)} we first observe that 
\begin{equation}\label{E12} 
d\gamma(U_1,U_2,L_1,L_2)=-\gamma([L_1,L_2],U_1,U_2)=0, 
\end{equation} 
where the first equality follows from $U_1,U_2\in \fz(\fl)$ and 
the second equality follows from $[L_1,L_2]\in\fl'$ and $\dim \fl'=2$. 
Combining now 
(\ref{E12}) with the cocycle condition for $(\alpha,\gamma)$ 
we obtain 
$$\langle\alpha(U_1,U_2),\alpha(L_1,L_2)\rangle + 
\langle\alpha(U_2,L_1),\alpha(U_1,L_2)\rangle + 
\langle\alpha(L_1,U_1),\alpha(U_2,L_2)\rangle=0.$$ 
Since {\it (i)} implies $\alpha(U_1,U_2)=0$ the first term vanishes and the 
assertion follows. 
\qed 

\begin{lm} \label{Beh2} 
Let $[\alpha,\gamma]\in\cH^2_Q(\fl,\fa)$ be 
admissible and choose $\alpha$ such that $\alpha(\fl,\fl)\subset\fa^\fl$. Let 
$\bar \fl$ be as in Lemma~\ref{L1}. If 
$K\in 
\fl'$ satisfies 
$\alpha(K,\bar \fl)=0$, then 
$\alpha(K,\fl)=0$. 
\end{lm} 
\proof Choose $X_1,X_2,X_3,Y,Z$ as in Lemma \ref{L1}. 
For $L_1,L_2\in\fl$ we have 
\begin{eqnarray*} 
\langle \alpha(K,L_1),\alpha(Y,L_2)\rangle &=& \langle 
\alpha(K,L_1),\alpha([X_1, 
X_2],L_2)\rangle\\ 
&=&\langle \alpha(K,L_1),\alpha([L_2,X_2],X_1)\rangle + \langle 
\alpha(K,L_1),\alpha([X_1,L_2],X_2)\rangle \\ 
&=&\langle \alpha(K,X_1),\alpha([L_2,X_2],L_1)\rangle + \langle 
\alpha(K,X_2),\alpha([X_1,L_2],L_1)\rangle\\ 
&=&0, 
\end{eqnarray*} 
where we first used the cocycle condition for $\alpha$ and then Lemma 
\ref{L3}. 
Similarly, we have 
$$\langle \alpha(K,L_1),\alpha(Z,L_2)\rangle = \langle 
\alpha(K,L_1),\alpha([X_1, 
X_3],L_2)\rangle=0.$$ 
This implies $\alpha(K,L_1)\perp\alpha(\fl',\fl)$. Now $(B_1)$ yields 
$\alpha(K,L_1)=0$ for 
all $L_1\in\fl$. 
\qed 

\begin{lm} \label{Beh3} 
Let $\fl$ be admissible and let $\bar\fl$ be as in Lemma \ref{L1}. Then 
$[L_1,L_2]=0$ 
holds for all 
$L_1,L_2\in\fl$ satisfying $[L_1,\bar \fl\,]=[L_2,\bar \fl\,]=0$. 
\end{lm} 
\proof  We choose a semi-simple orthogonal 
$\fl$-module $\fa$ such that there is an admissible cohomology class  
$[\alpha,\gamma]\in \cHQ$. We may assume $\alpha(\fl,\fl)\subset \fa^\fl$. From 
$d\alpha=0$ 
and $[\,\bar\fl,L_1]=[\,\bar\fl,L_2]=0$ we obtain 
$\alpha([L_1,L_2],X_{i})=0$. Lemma \ref{Beh2} 
now implies 
\begin{equation}\label{babe} 
\alpha([L_1,L_2],\cdot)=0. 
\end{equation} 

Take $X_1,X_2,X_3,Y,Z$ as in Lemma \ref{L1}. Using $[L_1,L_2]\in \fl'\subset 
\fz(\fl)$ we see that 
$$d\gamma([L_1,L_2],X_{1},X_{2},L_i) = 
-\gamma(Y,[L_1,L_2],L_i)$$ for $i=1,2.$ 
On the other hand (\ref{babe}) gives 
$$2 d\gamma([L_1,L_2],X_{1},X_{2},L_i) = \langle 
\alpha\wedge\alpha \rangle ([L_1,L_2],X_{1},X_{2},L_i) = 0$$ 
for $i=1,2.$ Hence we have 
\begin{equation}\label{a} 
\gamma(Y,[L_1,L_2],L_i)=0 
\end{equation} 
for $i=1,2$ and similarly we obtain 
\begin{equation}\label{b} 
\gamma(Z,[L_1,L_2],L_i)=0 
\end{equation} 
for $i=1,2$. Assume now that $[L_1,L_2]\not=0$. By 
(\ref{a}), 
(\ref{b}) and $[L_1,L_2]\in\Span\{Y,Z\}$ we obtain 
$\gamma(Y,Z,L_i)=0$ 
for $i=1,2$. This yields 
\begin{equation}\label{c} 
d\gamma(U,L_1,L_2,L)=-\gamma([L_1,L_2],U,L) 
\end{equation} 
for all $L\in\fl$ and $U\in \fl'$. On the other hand, Lemma \ref{L3}, {\it(ii)}, 
{\it(iii)} gives 
$\alpha(U,L_1)={\alpha(U,L_2)=0}$ and therefore 
\begin{equation}\label{d} 
d\gamma(U,L_1,L_2,L)=\langle\alpha(U,L),\alpha(L_1,L_2)\rangle. 
\end{equation} 
{}From (\ref{c}) and (\ref{d}) we get 
$$\gamma(L,[L_1,L_2],\cdot)=\langle 
\alpha(L_1,L_2),\alpha(L,\cdot)\rangle$$ 
as an element of $(\fl')^*$. 
Hence Condition $(A_{1})(ii)$ is satisfied for $L_0=[L_1,L_2]$, 
$A_0=-\alpha(L_1,L_2)$, 
$Z_0=0$. 
Since also $(A_{1})(i)$ holds by (\ref{babe}) and $[\alpha,\gamma]$ is 
admissible 
we get $[L_1,L_2]=0$, which is a contradiction. 
\qed 

\begin{lm} \label{Beh4} 
If $[\alpha,\gamma]\in\cHQ$ is admissible and $\alpha(\fl,\fl)\subset\fa^\fl$, 
then 
$\alpha([L,\fl],\fl)=0$ holds for all $L\in\fl$ satisfying $[L,\bar\fl\,]=0$. 
\end{lm} 
\proof 
Let $L\in\fl$ satisfy $[L,\bar\fl\,]=0$. By Lemma \ref{Beh2} it suffices to 
prove 
that $\alpha([L,L'],X)=0$ 
holds for all $L'\in\fl$ and $X\in\bar \fl$. Since $[X,L]=0$ the cocycle 
condition 
for $\alpha$ gives 
$$d\alpha(L,L',X)=-\alpha([L,L'],X)-\alpha([L',X],L)=0.$$ 
The assertion now follows since $\alpha([L',X],L)=0$ by Lemma 
\ref{L3} $(ii)$, $(iii)$. 
\qed 

\begin{lm}\label{L2} 
There exists a basis $X_1,X_2,X_3,\dots,X_l,Y,Z$ of $\fl$ such that 
\begin{eqnarray} 
&[X_1,X_2]=Y,\ [X_1,X_3]=Z,\ [X_2,X_3]=0&\label{E9}\\ 
& [X_1,X_4]=0,\ [X_2,X_4]=\lambda Z,\ \lambda\in\{0,1\},&\label{E10}\\ 
&[X_1,X_j]=[X_2,X_j]=0,\ j\ge 5.&\label{E11} 
\end{eqnarray} 
Let $X_1,X_2,X_3,\dots,X_l,Y,Z$ be such a basis and let $j_0$ be such that 
$$[X_1,X_j]=[X_2,X_j]=[X_3,X_j]=0$$ for all $j\ge j _0$. Then we have 
$[X_r,X_s]=0$ for all $r,s\ge j_0$. 
\end{lm} 
\proof 
We choose $\bar \fl$, a basis $X_1,X_2,X_3$ of $\bar\fl$ and a basis $Y,Z$ of 
$\fl'$ as in Lemma \ref{L1}. 
Because of $[X_1,X_2]=Y$ and $[X_1,X_3]=Z$ one can find a complementary 
vector space $W$ 
of $\Span\{X_1,X_2,X_3,Y,Z\}$ in $\fl$ such that $[X_1,W]=0$ and 
$[X_2,W]\subset \RR\cdot Z$ (choose a basis of an arbitrary complement and 
change each 
basis vector by a suitable linear combination of $X_1$, $X_2$ and $X_3$). If 
$[X_2,W]=0$, 
then we can choose an arbitrary basis $X_4,\dots,X_l$ of $W$ and (\ref{E10}) 
and 
(\ref{E11}) are satisfied for $\lambda=0$. If $[X_2,W]\not=0$ then we can 
choose a basis 
$X_4,\dots,X_l$ of $W$ such that $[X_2,X_4]=Z$ and $[X_2,X_j]=0$ for 
$j=1,\dots,l$. 

The second statement follows from Lemma \ref{Beh3}.
\qed 

Now we fix a basis of $\fl$ which satisfies the conditions of Lemma 
\ref{L2}.

\subsubsection*{Case II.1: $\lambda=0$} 
Let $X_1,X_2,X_3,\dots,X_l,Y,Z$ be a basis of $\fl$ satisfying 
(\ref{E9}), (\ref{E10}) and (\ref{E11}) with $\lambda=0$. We define 
$W:=\Span 
\{X_{4},\dots,X_{l}\}$. 

\subsubsection*{Case II.1.1: $[X_{3},W]=0$} 
{\bf Claim.} An admissible Lie algebra $\fl$ which satisfies the conditions 
of case II.1.1 
is isomorphic to $\fg_{5,2}\oplus \RR^k$, \\[2ex] 
\proof In this case we have $[X_{1},W]=[X_{2},W]=[X_{3},W]=0$ by 
assumption and $[W,W]=0$ by Lemma \ref{Beh3}. 
\qed 

\subsubsection*{Case II.1.2: $\dim [X_{3},W]=1$}

{\bf Claim.} An admissible Lie algebra $\fl$ which satisfies the conditions 
of case II.1.2 is isomorphic to $\fg_{6,4}\oplus \RR^k$ or to 
$\fh(1)\oplus\fh(1)$.\\[2ex] 
\proof We may assume $[X_{3},X_{4}]=cY+dZ\not=0$, $c,d\in\RR$ and 
$[X_{3},X_{r}]=0$ for $r\ge 5$. Let us first consider the case $d\not=0$. 
Replacing 
$X_{3}$, $X_{4}$, and $Z$ by 
$$X_{3}':=dX_{3}+cX_{2},\ X_{4}':=1/d\cdot X_{4},\mbox{ and } 
Z':=cY+dZ,$$ 
respectively, we see that we may assume $c=0$ and $d=1$. Hence we have 
a basis $X_1,X_2,X_3,\dots,X_l,Y,Z$ of $\fl$ satisfying (\ref{E9}), 
(\ref{E10}), (\ref{E11}), $[X_{3},X_{4}]=Z$ and $[X_{3},X_{r}]=0$ for 
$r\ge 5$. 
We will prove that the admissibility of $\fl$ implies 
$[X_{4},X_{r}]=0$ for $r\ge 5$. Assume first that there is a vector $X\in 
\Span\{X_5,\dots,X_l\}$ such that $[X_{4},X]=aY+bZ$, $a\not=0$ and $b\not=0$. 
Let 
$[\alpha,\gamma]\in\cHQ$ be 
admissible and choose $\alpha$ such that $\alpha(\fl,\fl)\subset \fa^\fl$. Then 
the cocycle condition for 
$\alpha$ 
yields 
\begin{eqnarray*} 
d\alpha(X_{1},X_{2},X_{3}) &=& -\alpha(Y,X_{3})+\alpha(Z,X_{2})\ =\ 0  
\\ 
d\alpha(X_{1},X_{2},X_{4}) &=& -\alpha(Y,X_{4})\ =\ 0 \\ 
d\alpha(X_{1},X_{3},X_{4}) &=& -\alpha(Z,X_{4})- \alpha(Z,X_{1})\ =\ 0  
\\ 
d\alpha(X_{2},X_{3},X_{4}) &=& -\alpha(Z,X_{2})\ =\ 0. 
\end{eqnarray*} 
Moreover, Lemma \ref{Beh4} for $\bar \fl=\Span\{X_1,X_2,X_3\}$, $L=X$ yields 
$\alpha(a Y +bZ,\cdot)=0$. Since 
$a\not=0$ and $b\not=0$, this equation together with the cocycle 
conditions above implies $\alpha(\fl',X_i)=0$ for $i\le 4$ and Lemma~\ref{L3}, 
{\it(ii)}, {\it(iii)} now gives $\alpha(\fl',\fl)=0$. In particular we obtain 
$d\gamma(\fl',\fl,\fl,\fl)=0$. From this condition we obtain 
\begin{eqnarray*} 
d\gamma(Z,X_{1},X_{2},X_{j}) &=& - \gamma(Y,Z,X_{j})\ =\ 0,\ j\ge 3 \\ 
d\gamma(Y,X_{1},X_{3},X_{4}) &=& \gamma(Y,Z,X_{4}) + 
        \gamma(Y,Z,X_{1}) \ =\ 0 \\ 
d\gamma(Y,X_{2},X_{3},X_{4}) &=& \gamma(Y,Z,X_{2}) \ =\ 0, 
\end{eqnarray*} 
hence $ \gamma(Y,Z,\cdot)=0$. But then $[\alpha,\gamma]$ does not 
satisfy Condition $(A_{1})$. This is a contradiction to the 
admissibility of $[\alpha,\gamma]$. Hence $[X_4,X_r]\in \RR Z$ for all $r\ge 5$ 
or 
$[X_4,X_r]\in \RR Y$ for all $r\ge 5$. If we are in the first case and if there 
is an 
$s\ge 5$ such that $[X_4,X_s]\not=0$, then we may assume $[X_4,X_5]=Z$. If we 
define 
$\bar\fl =\Span\{X_1,X_2,X_{3}+X_{5}\}$, then $\bar\fl$, $L_1=X_4$ and 
$L_2=X_5$ 
satisfy the assumptions of Lemma \ref{Beh3}. But $[X_{4},X_{5}]=Z\not=0$ yields 
a 
contradiction. Similarly, if $[X_4,X_r]\in \RR Y$ for all $r\ge 5$ and if there 
is an 
$s\ge 5$ such that $[X_4,X_s]\not=0$, then we may assume $[X_4,X_5]=Y$. Consider 
now 
$\bar\fl:=\Span\{X_{4},X_{3},X_{5}-X_{2}\}$. Then $\bar\fl$, $L_1=X_1+X_4$ and 
$L_2=X_2$ satisfy the assumptions of Lemma \ref{Beh3}, but 
$[X_1+X_4,X_2]=Y\not=0$, a 
contradiction. 
We deduce $[X_{4},X_{r}]=0$ for $r\ge 5$. We conclude that in 
case $d\not=0$ the Lie algebra $\fl$ is isomorphic to 
$$\{ [X_{1},X_{2}]=Y,\ [X_{1},X_{3}]=Z,\ [X_{3},X_{4}]=Z\} \oplus 
\RR^{k}.$$ Putting $X_{1}':=X_{1}+ X_{4}$ we see 
$$\fl\cong \{ [X_{1}',X_{2}]=Y,\ [X_{3},X_{4}]=Z\} \oplus 
\RR^{k}\cong \fh(1)\oplus \fh(1)\oplus 
\RR^{k}.$$ 

Now we consider the case $d=0$. We may assume $c=1$. 
Now we have 
a basis $X_1,X_2,X_3,\dots,X_l,Y,Z$ of $\fl$ satisfying (\ref{E9}), 
(\ref{E10}), (\ref{E11}), $[X_{3},X_{4}]=Y$ and $[X_{3},X_{r}]=0$ for 
$r\ge 5$. We will prove that $[X_{4},X_{r}]=0$ holds for $r\ge5$. Assume that 
this is 
not true. Then we have without loss of generality $[X_{4},X_{5}]=aY+bZ\not=0$. 
If 
$b\not=0$, then we replace 
$X_{3}$ by $X_{3}':=X_{3}+X_{5}$. Then the basis 
$X_1,X_2,X_3',X_{4},\dots,X_l,Y,Z$ 
satisfies (\ref{E9}), 
(\ref{E10}), (\ref{E11}) with $\lambda=0$, $[X_{3}',X_{4}]=(1-a)Y 
-bZ$, and $[X_{3}',X_{r}]=0$ for $r\ge 5$. Thus we are in the above 
case where $d\not=0$. This implies $[X_{4},X_{5}]=0$, a 
contradiction. If $b=0$, then we may assume $a=1$. Again $\bar \fl=\Span 
\{X_1,X_2,X_{3}+X_{5}\}$, $L_1=X_{4}$ and $L_2=X_{5}$ satisfy the assumptions 
of 
Lemma \ref{L2}, but $[X_4,X_5]\not=0$, a 
contradiction. Therefore we have $[X_{4},X_{r}]=0$ for $r\ge5$, thus $\fl$ is 
isomorphic to 
$$\{ [X_{1},X_{2}]=Y,\ [X_{1},X_{3}]=Z,\ [X_{3},X_{4}]=Y\} \oplus 
\RR^{k} \cong \fg_{6,4}\oplus 
\RR^{k}.$$ 
\qed 
\subsubsection*{Case II.1.3: $\dim [X_{3},W]=2$}

{\bf Claim.} 
A Lie algebra $\fl$ which satisfies the conditions of case II.1.3 is not 
admissible. 
\\[2ex] 
\proof 
Obviously we may assume that $X_1,\dots,X_l,Y,Z$ is a basis of $\fl$ which 
satisfies (\ref{E9}), 
(\ref{E10}), (\ref{E11}) with $\lambda=0$, $[X_{3},X_{4}]=Y$, $[X_{3},X_{5}]=Z$, 
and $[X_{3},X_{r}]=0$ for $r>5$. Moreover, we have $[X_4,X_5]=yY+zZ$ for 
suitable $y,z\in\RR$. Assume first that $z\not=0$. Let 
$[\alpha,\gamma]\in\cH^2_Q(\fl,\fa)$ be admissible and choose $\alpha$ such that 
$\alpha(\fl,\fl)\subset \fa^\fl$. 
Then we have 
\begin{eqnarray} 
d\alpha(X_{1},X_{2},X_{3}) &=& -\alpha(Y,X_{3})+\alpha(Z,X_{2})\ =\ 0  
\label{q1}\\ 
d\alpha(X_{2},X_{3},X_{4}) &=& -\alpha(Y,X_{2})\ =\ 0 \label{q2} \\ 
d\alpha(X_{2},X_{3},X_{5}) &=& -\alpha(Z,X_{2})\ =\ 0 \label{q3} \\ 
d\alpha(X_{1},X_{3},X_{5}) &=& -\alpha(Z,X_{5})-\alpha(Z,X_{1})\ =\ 
0\label{q4}\\ 
d\alpha(X_{1},X_{3},X_{4}) &=& -\alpha(Z,X_{4})-\alpha(Y,X_{1})\ =\ 0.\label{q5} 
\end{eqnarray} 
By Lemma \ref{L3} {\it (ii)} we have $\alpha(Y,X_j)=0$ for $j\ge 4$. Together 
with 
(\ref{q1}) -- (\ref{q3}) this yields $\alpha(Y,X_j)=0$ for $j\ge 2$. 
Lemma \ref{L3} {\it (iv)} gives 
$$\langle\alpha(Y,X_1),\alpha(Z,X_j)\rangle \ =\ \langle 
\alpha(Y,X_j),\alpha(Z,X_1)\rangle \ =\ 0$$ 
for $j\ge 2$ and therefore also 
\begin{eqnarray*} 
\langle \alpha(Y,X_1),\alpha(Z,X_1)\rangle &= &-\langle 
\alpha(Y,X_1),\alpha(Z,X_5)\rangle \ =\ 0\\ 
\langle \alpha(Y,X_1),\alpha(Y,X_1)\rangle &= &-\langle 
\alpha(Y,X_1),\alpha(Z,X_4)\rangle \ =\ 0 
\end{eqnarray*} 
where we used (\ref{q4}) and (\ref{q5}). We obtain  
$\alpha(Y,X_1)\perp\alpha(\fl',\fl)$. Now we use that the admissibility 
condition 
$(B_1)$ implies that $\alpha(\fl,\fl')$ 
is non-degenerate. Hence $\alpha(Y,X_1)=0$ and, consequently, $\alpha(Y,\fl)=0$. 
Now we get $\alpha(Z,X_2)=\alpha(Z,X_4)=0$ from (\ref{q3}) and (\ref{q5}). 
Moreover, 
$$d\alpha(X_1,X_4,X_5)\ =\ \alpha(X_1,yY+zZ)\ =\ 0$$ 
gives $\alpha(Z,X_1)=0$ because of $z\not=0$. Thus $\alpha(Z,X_5)=0$ by 
(\ref{q4}). 
Finally, 
$$d\alpha(X_3,X_4,X_5)\ =\ -\alpha(Y,X_5)-\alpha(yY+zZ,X_3)+\alpha(Z,X_4)\ =\ 
0$$ 
yields $\alpha(Z,X_3)=0$. Since, moreover, $\alpha(Z,X_k)=0$ for $k\ge6$ by 
Lemma 
\ref{L3} {\it (ii)} we obtain also $\alpha(Z,\fl)=0$. 

Because of $2d\gamma=\langle \alpha\wedge\alpha\rangle$ we now get 
\begin{eqnarray*} 
d\gamma(Z,X_{1},X_{3},X_{4}) &=& - \gamma(Y,Z,X_{1})\ =\ 0, \\ 
d\gamma(Y,X_{1},X_{2},X_{3}) &=& -\gamma(Y,Z,X_{2}) \ =\ 0 ,\\ 
d\gamma(Z,X_{1},X_{2},X_{3}) &=& -\gamma(Y,Z,X_{3}) \ =\ 0,\\ 
d\gamma(Z,X_{1},X_{2},X_{k}) &=& -\gamma(Y,Z,X_{k}) \ =\ 0,\ k\ge 4. 
\end{eqnarray*} 
Consequently, $\gamma(Y,Z,\fl)=0$. This together with $\alpha(\fl',\fl)=0$ 
yields a 
contradiction to admissibility. Hence $z=0$ and $[X_4,X_5]=yY$. However, now we 
can 
apply Lemma \ref{Beh3} to $\bar\fl:=\Span\{X_3,X_4+yX_2,X_5\}$, $L_1=X_1+X_5$, 
$L_2=X_2$ and we obtain a contradiction to $[X_1+X_5,X_2]=Y\not=0$. 
\qed

\subsubsection*{Case II.2: $\lambda=1$} 

{\bf Claim.} An admissible Lie algebra $\fl$ which satisfies the conditions 
of case II.2 and which does not have a basis satisfying already the conditions 
of case II.1 is isomorphic to 
$\fg_{6,5}\oplus \RR^k$.\\[2ex] 
Let $X_1,\dots,X_l,Y,Z$ be a basis of $\fl$ which satisfies (\ref{E9}), 
(\ref{E10}), and (\ref{E11}) with $\lambda=1$. 
\begin{lm} \label{ix} 
If $[\alpha,\gamma]\in\cHQ_\sharp$ and $\alpha(\fl,\fl)\subset \fa^\fl$, 
then 
$\alpha([X_3,X_j],\cdot)=0$ for 
all $j\ge 5$. 
\end{lm} 
\proof 
By Lemma \ref{Beh2} it suffices to prove $\alpha([X_3,X_j],X_i)=0$ for all $j\ge 
5$ and $i=1,2,3$. For $i=2$ this follows obviously from the cocycle condition 
for $\alpha$. 

Now we consider $i=3$. Using $[X_1,X_2]=Y$, the cocycle condition for $\alpha$ 
and Lemma \ref{L3} {\it (iv)} we see that 
\begin{eqnarray*} 
\langle \alpha([X_3,X_j],X_3),\alpha(Y,L)\rangle &=& 
\alpha([X_3,X_j],L),\alpha(Y,X_3)\rangle\\ 
&=& \langle \alpha([X_3,X_j],L),\alpha([X_1,X_2],X_3)\rangle\\ 
&=& \langle \alpha([X_3,X_j],L),\alpha([X_1,X_3],X_2)\rangle\\ 
&=& \langle \alpha([X_3,X_j],X_2),\alpha([X_1,X_3],L)\rangle \ =\ 0, 
\end{eqnarray*} 
where the last equality follows from the above considerations for $i=2$. 
Similarly (using now $[X_2,X_4]=Z$) we obtain 
\begin{eqnarray*} 
\langle \alpha([X_3,X_j],X_3),\alpha(Z,L)\rangle &=& 
\alpha([X_3,X_j],L),\alpha(Z,X_3)\rangle\\ 
&= &\langle \alpha([X_3,X_j],L),\alpha([X_2,X_4],X_3)\rangle\\ 
&= &\langle \alpha([X_3,X_j],L),\alpha([X_3,X_4],X_2)\rangle\\ 
&= &\langle \alpha([X_3,X_j],X_2),\alpha([X_3,X_4],L)\rangle \ =\ 0. 
\end{eqnarray*} 
Now we use that the admissibility condition $(B_1)$ implies that 
$\alpha(\fl,\fl')$ 
is non-degenerate. This gives $\alpha([X_3,X_j],X_3)=0$. 

Finally we consider the case $i=1$. Note first that the cocycle condition for 
$\alpha$ implies 
\begin{equation} 
\label{oben1} 
\alpha([X_3,X_j],X_1)=-\alpha([X_1,X_3],X_j)+\alpha([X_1,X_j],X_3)= 
-\alpha(Z,X_j)\,. 
\end{equation} 
Using now Lemma~\ref{L3}, $(ii)$ and $(iv)$ we obtain 
$$\langle \alpha([X_3,X_j],X_1), \alpha(Y,L)\rangle = 
-\langle \alpha(Z,X_j), \alpha(Y,L)\rangle = 
-\langle \alpha(Z,L), \alpha(Y,X_j)\rangle = 0$$ 
for all $L\in\fl$. Consequently, $\alpha([X_3,X_j],X_1)\perp\alpha(Y,\fl)$. 
Now we will prove that also $\alpha([X_3,X_j],X_1)\perp\alpha(Z,\fl)$ and thus 
$\alpha([X_3,X_j],X_1)\perp\alpha(\fl',\fl)$ holds. By  
$(B_1)$ this will give $\alpha([X_3,X_j],X_1)=0$. We observe that 
\begin{eqnarray} 
\lefteqn{\hspace{-1.5cm}\langle \alpha([X_3,X_j],X_1),\alpha(Z,L)\rangle = 
-\langle \alpha(Z,X_j),\alpha(Z,L)\rangle \,=\, -\langle 
\alpha(Z,X_j),\alpha([X_2,X_4],L)\rangle }\nonumber\\ 
&=&-\langle \alpha(Z,X_j),\alpha([X_2,L],X_4)\rangle -\langle 
\alpha(Z,X_j),\alpha([L,X_4],X_2)\rangle, \label{oben2} 
\end{eqnarray} 
where we used (\ref{oben1}) and the cocycle condition for $\alpha$. By Equation 
(\ref{oben1}) and 
Lemma~\ref{L3}~$(iv)$ we have 
\begin{eqnarray*} 
-\langle \alpha(Z,X_j),\alpha([L,X_4],X_2)\rangle&=&\langle 
\alpha([X_3,X_j],X_1),\alpha([L,X_4],X_2)\rangle\\ 
&=&\langle \alpha([X_3,X_j],X_2),\alpha([L,X_4],X_1)\rangle\ =\ 0 
\end{eqnarray*} 
since $\alpha([X_3,X_j],X_2)=0$. Hence the last term in (\ref{oben2}) vanishes 
and we get 
$$ 
\langle \alpha([X_3,X_j],X_1),\alpha(Z,L)\rangle = -\langle 
\alpha(Z,X_j),\alpha([X_2,L],X_4)\rangle 
= c\langle\alpha(Z,X_j),\alpha(Z,X_4)\rangle 
$$ 
for some real number $c\in\RR$ since $[X_2,L]\in\Span\{Y,Z\}$ and since 
$$\langle\alpha(Z,X_j),\alpha(Y,X_4)\rangle 
=\langle\alpha(Z,X_4),\alpha(Y,X_j)\rangle=0 $$ 
by Lemma~\ref{L3}~{\it (ii)} and $(iv)$. Furthermore, we have 
$$ 
\langle 
\alpha(Z,X_j),\alpha(Z,X_4)\rangle = \langle 
\alpha(Z,X_j),\alpha([X_1,X_3],X_4)\rangle\\ 
=- \langle \alpha(Z,X_j),\alpha([X_3,X_4],X_1)\rangle\,. 
$$ 
Since we already know that 
$\alpha(Z,X_j)\perp\alpha(Y,X_1)$ the last equation implies that in order to 
prove 
$\langle \alpha([X_3,X_j],X_1),\alpha(Z,L)\rangle=0$ it suffices to show 
$\langle \alpha(Z,X_j),\alpha(Z,X_1)\rangle=0$. However, this follows from 
Lemma~\ref{L3}, $(ii)$ and $(iv)$: 
\begin{eqnarray*} 
\langle \alpha(Z,X_j),\alpha(Z,X_1)\rangle&=& \langle 
\alpha(Z,X_j),\alpha([X_2,X_4],X_1)\rangle \, =\, -\langle 
\alpha(Z,X_j),\alpha([X_1,X_2],X_4)\rangle \\ &=& 
-\langle\alpha(Z,X_j),\alpha(Y,X_4)\rangle\ =\ 
-\langle\alpha(Y,X_j),\alpha(Z,X_4)\rangle \ =\ 0\,. 
\end{eqnarray*} 
\qed 

\begin{lm}\label{Lx} 
If $\fl$ is admissible, then $[X_3,W']=0$ for $W':=\Span\{X_5,\dots,X_l\}$. 
\end{lm} 
\proof 
Suppose $[X_3,W']\not=0$. Then we may assume that besides 
(\ref{E9}), (\ref{E10}), and (\ref{E11}) with $\lambda=1$ our basis satisfies 
also 
\begin{equation}\label{E35} 
[X_3,X_5]=uY+vZ\not=0\,. 
\end{equation} 
Take $[\alpha,\gamma]\in\cHQ_\sharp$ such that $\alpha(\fl,\fl)\subset \fa^\fl$. 
From Lemma \ref{ix} we know that 
\begin{equation}\label{xx} 
\alpha(uY+vZ,\cdot)=0\,, 
\end{equation} 
which implies 
\begin{eqnarray} 
0\ =\ d\gamma(uY+vZ,X_1,X_2,X_3) &=& 
-\gamma(Y,uY+vZ,X_3)+\gamma(Z,uY+vZ,X_2)\quad \label{x1}\\ 
0\ =\ d\gamma(uY+vZ,X_1,X_2,X_4) &=& 
-\gamma(Y,uY+vZ,X_4)-\gamma(Z,uY+vZ,X_1)\quad \label{x2}\\ 
0\ =\ d\gamma(uY+vZ,X_1,X_2,X_k) &=& -\gamma(Y,uY+vZ,X_k),\ k\ge5 \label{x3}\,. 
\end{eqnarray} 
Furthermore, we have 
\begin{eqnarray*} 
d\alpha(X_1,X_2,X_3)&=&-\alpha(Y,X_3)+\alpha(Z,X_2)\ =\ 0\\ 
d\alpha(X_1,X_2,X_4)&=&-\alpha(Y,X_4)-\alpha(Z,X_1)\ =\ 0\,. 
\end{eqnarray*} 
Let us first consider the case $v\not=0$ in (\ref{E35}). Replacing 
$X_3,X_4,X_5$ and $Z$ by 
$$X_3':= vX_3+uX_2 ,\ X_4':= vX_4-uX_1,\ X_5':=(1/v)\cdot X_5,\ Z':=uY+vZ $$ 
we see that we may assume $u=0$ and $v=1$ in (\ref{E35}), i.e.~$[X_3,X_5]=Z$. 
Then 
(\ref{xx}) says $\alpha(Z,\cdot)=0$, hence Lemma~\ref{L3}~$(ii)$ and the above 
equations for $\alpha$ imply $\alpha(Y,X_j)=0$ for $j\ge3$. 
Equations 
(\ref{x1}), (\ref{x2}), and 
(\ref{x3}) imply $\gamma(Y,Z,X_j)=0$ for $j\ge3$. Using all this we obtain 
$$ 
\gamma(Y,Z,X_i)=d\gamma(Y,X_{i},X_{3},X_{5}) = \langle 
\alpha(Y,X_{i}),\alpha(X_{3},X_{5})\rangle,\quad i=1,2. 
$$  
In particular, the data $L_0=Z$, $A_0=-\alpha(X_3,X_5)$, and $Z_0=0$ satisfy the 
conditions $(i)$ and $(ii)$ of $(A_{1})$. Hence $Z=0$ by 
admissibility, which is a contradiction. 

If $v=0$, then we may assume $u=1$, i.e.~$[X_3,X_5]=Y$. Then (\ref{xx}) implies 
$\alpha(Y,\cdot)=0$. The above equations for $\alpha$ now give 
$\alpha(Z,X_j)=0$ for $j=1,2$. Hence $\alpha(X_1,\fl)=\alpha(X_2,\fl)=0$ and 
therefore 
\begin{eqnarray*} 
&0=d\alpha(X_2,X_4,X_3)=-\alpha(Z,X_3)\;&\\ 
&0=d\alpha(X_1,X_3,X_k)=-\alpha(Z,X_k),&k\ge4. 
\end{eqnarray*} 
This implies $\alpha(\fl',\fl)=0$. 
{}From (\ref{x1}) and (\ref{x2}) we obtain $\gamma(Y,Z,X_j)=0$ for $j=1,2$. 
Using this we get 
$$ 
2\gamma(Y,Z,X_j)=2d\gamma(Z,X_{1},X_{2},X_{j}) = \langle 
\alpha\wedge\alpha\rangle (Z,X_{1},X_{2},X_{j})=0 
$$for all $j\ge3$. 
Hence $\gamma(\fl',\fl',\fl)=0$. Again we obtain a 
contradiction to the admissibility condition $(A_{1})$. 
\qed 

\begin{lm} \label{Lxx} 
If $\fl$ is admissible, then $[X_4,W']=0$.  
\end{lm} 
\proof 
Recall that $X_{1},\ldots,X_l,Y,Z$ is a basis of $\fl$ which satisfies 
(\ref{E9}), (\ref{E10}), (\ref{E11}) with $\lambda=1$. Therefore the 
basis $X_{1}',\ldots,X_{l}',Y',Z'$ of 
$\fl$ defined by 
$$X_{1}':=X_{2},\ X_{2}':=X_{1},\ X_{3}':=X_{4},\ X_{4}':=X_{3},\ 
X_{j}':=X_{j},\ j\ge5,\ Y':=-Y,\ Z':=Z$$ 
also satisfies (\ref{E9}), (\ref{E10}), (\ref{E11}) with $\lambda=1$. 
Now Lemma \ref{Lx} says that 
$$[X_4,X_{j}]=[X_{3}',X_{j}']=0,\ j\ge5.$$ 
\qed 

{\sl Proof of the Claim.} We know from Equations (\ref{E9}), (\ref{E10}), 
(\ref{E11}), Lemma \ref{Lx} and Lemma~\ref{Lxx} that $\fl$ is isomorphic 
to $\fl_{1}\oplus \RR^{k}$, where 
$$\fl_{1}=\{[X_{1},X_{2}]=Y,\ [X_{1},X_{3}]=Z,\ [X_{2},X_{4}]=Z,\ 
[X_{3},X_{4}]=aY+bZ\} $$ 
for suitable $a,b\in\RR$.

Assume $a=0$. Then the basis 
$$X_{1},\ 
X_{2}':=X_{3}-bX_{2}, \ X_{3}':=X_{2},\ X_{4},\ Y':=Z-bY,\ Z':=Y$$ 
of $\fl_1$ 
together with a basis of $\RR^{k}$ satisfies the conditions of case 
II.1 which contradicts our assumption on $\fl$. 
Hence $a\not=0$. 

Replacing $X_{2},X_{4},Y$ by 
$$X_{2}':=X_{2}+ 
(b/2a)\cdot X_{3},\ X_{4}':=X_{4} 
+(b/2)\cdot X_{1},\ Y':=Y+(b/2a)\cdot Z$$ 
we obtain a basis $X_{1},X_{2}',X_{3},X_{4}',Y,Z$ of $\fl_{1}$ 
satisfying 
$$[X_{1},X_{2}']=Y',\ [X_{1},X_3]=Z,\ [X_{1},X_{4}']=[X_{2}',X_3]=0,\ 
[X_{2}',X_{4}']=\lambda'Z,\ [X_3,X_{4}']=aY'$$ 
where $\lambda'=1+b^{2}/4a$. Since by assumption $\fl$ does not have 
a basis satisfying the conditions of case II.1 we have 
$\lambda'\not=0$. Hence, we may obviously assume $\lambda'=1$. Putting 
$\mu=\sqrt{|a|}$ and 
$$X_{1}'=\mu X_{1},\ X_{3}'=(1/\mu)\cdot X_{3},\ \bar Y=\mu Y'$$ 
we obtain a basis $X_{1}',\ldots,X_{4}',\bar Y,Z$ of $\fl_{1}$ which 
satisfies 
$$[X_{1}',X_{2}']=\bar Y,\ [X_{1}',X_3']=Z,\ [X_{1}',X_{4}']=[X_{2}',X_3']=0,\ 
[X_{2}',X_{4}']=Z,\ [X_3',X_{4}']=\epsilon \bar Y$$ 
with $\epsilon=\pm1$. If $\epsilon=-1$, then $\fl_1\cong\fg_{6,5}$. If $\epsilon 
=1$, then 
$$\fl_1\cong\{[\tilde X_1,\tilde X_2]=\tilde Y,\  [\tilde X_3,\tilde X_4]=\tilde 
Z\}\cong\fh(1)\oplus\fh(1)$$
for 
\begin{eqnarray*}&\tilde X_1=\textstyle{\frac 12}(X_1'-X_4'),\ \tilde X_2=\frac 
12(X_2'+X_3'),\ \tilde X_3=\frac 12(X_1'+X_4'),\ \tilde X_4=\frac 
12(X_2'-X_3'),&\\
&\tilde Y=\frac 12(\bar Y+Z),\ \tilde Z=\frac 12(\bar Y-Z)\,.&
\end{eqnarray*}
\qed 

\section{Nilpotent metric Lie algebras of dimension $\le 10$} 
\label{S4}
Recall that a metric Lie algebra is called indecomposable if it is not the 
direct sum of two non-trivial metric Lie algebras (see also 
\cite{KO1}). 
In this section we will determine all indecomposable nilpotent metric Lie 
algebras of dimension $\le 10$ (up to isomorphisms). 

Let us first consider the following construction. Let $\fl$ be a 
nilpotent Lie algebra and let $(\fa,\ip_{\fa})$ be a pseudo-Euclidean 
vector space which we consider as a trivial orthogonal $\fl$-module. 
Let $\fd$ be the vector space $\fl^{*}\oplus\fa\oplus\fl$. 
Take $(\alpha,\gamma)\in\cZQ$ and define a bilinear map $\lb:\fd 
\times \fd\rightarrow \fd$ by
$$\begin{array}{l}
    [\fl^{*}\oplus\fa,\fl^{*}\oplus\fa ]=0,\ [\fl,\fa]=0   \\[1ex]
    [L,Z]=\ad^{*}(L)(Z)   \\[1ex]
    [A,L]=\langle A,\alpha(L,\cdot)\rangle \\[1ex]
    [L_{1},L_{2}]=\gamma 
    (L_{1},L_{2},\cdot)+\alpha(L_{1},L_{2})+[L_{1}L_{2}]_{\fl} 
\end{array}$$
for all $L,L_{1},L_{2}\in\fl$, $A\in\fa$, and $Z\in
\fl^{*}$. Moreover we define an inner product $\ip$ on $\fd$ by 
\begin{eqnarray*}
 \langle Z_1+A_1+L_1,Z_2+A_2+L_2\rangle&:=& \langle A_1,A_2\rangle_\fa
+Z_1(L_2) +Z_2(L_1) 
\end{eqnarray*}
for $Z_1,\,Z_2\in \fl^*$, $A_1,\,A_2\in \fa$ and
$L_1,\,L_2\in \fl$. Then it is not hard to prove that 
$\dd:=(\fd,\lb,\ip)$ 
is a nilpotent metric Lie algebra
(see also \cite{KO2} for the case of a general metric Lie algebra).

Let $\fl_{i}$, $i=1,2$ be Lie algebras and let $\fa_{i}$, $i=1,2$ be 
pseudo-Euclidean vector spaces which we consider as trivial 
orthogonal $\fl_{i}$-modules. Consider a pair $(S,U)$ consisting of 
a homomorphism 
$S:\fl_{1}\rightarrow \fl_{2}$ and an isometry $U:\fa_{2}\rightarrow 
\fa_{1}$. Then $(S,U)^{*}: C^{p}(\fl_{2},\fa_{2})\rightarrow 
C^{p}(\fl_{1},\fa_{1})$ induces a map 
$(S,U)^{*}:\cH_{Q}^{p}(\fl_{2},\fa_{2})\rightarrow 
\cH_{Q}^{p}(\fl_{1},\fa_{1})$.

In particular, $G_{(\fl,\fa)}:= \Aut(\fl)\times O(\fa,\ip_\fa)$ acts on $\cHQ$.

\begin{de} Let $\fl$ be a 
nilpotent Lie algebra and let $(\fa,\ip_{\fa})$ be a pseudo-Euclidean 
vector space considered as a trivial $\fl$-module. A cohomology class 
$\ph\in\cHQ$ is called decomposable if there are decompositions 
$\fa=\fa_1\oplus\fa_2$ and $\fl=\fl_1\oplus\fl_2$, at least one of them being 
non-trivial and cohomology classes $\ph_i\in\cH_Q^2(\fl_i,\fa_i)$, $i=1,2$ such 
that $\ph=(q_1,j_1)^*\ph_1+(q_2,j_2)^*\ph_2$, where $q_i:\fl\rightarrow \fl_i$ 
are the projections and $j_i:\fa_i\rightarrow \fa$ are the inclusions. Here we 
consider  $\fa_i$ as trivial $\fl_i$-modules.
We denote the subset of indecomposable admissible cohomology classes in $\cHQ$ 
by $\cHQ_0$.
\end{de}
One can check easily that $\cHQ_0$ is invariant with respect to the action of 
$G_{(\fl,\fa)}$ on $\cHQ$. The classification scheme (\ref{x}) now gives

\begin{pr} \label{P3} The set of isomorphism classes of nilpotent metric Lie 
algebras of dimension at most 10 is in bijective correspondence with 
$$\bigcup_{\fl\in\frak L}\bigcup_{\fa\in\frak A_\fl}\cHQ_0/G_{(\fl,\fa)},$$ 
where 
$\frak L$ is the set of isomorphism classes of nilpotent Lie algebras of 
dimension at most 5 and for a fixed $\fl\in\frak L$ 
the set $\frak A_\fl$ consists of all isometry classes of pseudo-Euclidean 
vector spaces of dimension at most $10-2\dim\fl$ which we consider as 
equivalence classes of trivial orthogonal $\fl$-modules.
\end{pr}

In the following we will often abbreviate $G_{(\fl,\fa)}$ to $G$. Furthermore, 
we will use the following conventions. An orthonormal basis  of a 
pseudo-Euclidean vector space $(\fa,\ip_\fa)$ is a basis $A_1,\dots,A_{p+q}$ 
consisting of pairwise orthogonal vectors satisfying $\langle 
A_i,A_i\rangle_\fa=-1$ for $1\le i\le p$ and $\langle A_i,A_i\rangle_\fa=1$ for 
$p+1\le i\le p+q$. The pair $(p,q)$ is called signature of $\fa$. We denote the 
standard pseudo-Euclidean vector space of signature $(p,q)$ by $\RR^{p,q}$. A 
Witt basis of $\RR^{1,1}$ is a basis $A_1,A_2$, where $A_1,A_2$ are isotropic 
and $\langle A_1,A_2\rangle =1$. 
\begin{pr} \label{P4}
If $\fl$ is nilpotent and if $\dim \fl=5$ and $\fa=0$, then 
$\cHQ_0\not=\emptyset$ 
implies 
$\fl=\RR^5$ or $\fl=\fg_{5,2}$. 
\end{pr} 
\proof 
Let $[0,\gamma]\in\cHQ_0$ be such that $\gamma\not=0$. Then we know from $(A_k)$ 
that 
$\dim 
\fl^{k+1}\not=1$ holds for all $k\ge 0$. Since $\fl$ is nilpotent the 
codimension of 
$\fl^2$ 
in $\fl$ cannot be 1, since otherwise $\fl^3=[\fl,\fl^2]=[\fl,\fl]=\fl^2$ yields 
a 
contradiction. Hence we have only the following possibilities: 
\begin{itemize} 
\item[(i)] $\dim \fl^2=0$, \vspace{-1ex} 
\item[(ii)] $\dim \fl^2=2, \dim \fl^3=0$,\vspace{-1ex} 
\item[(iii)] $\dim \fl^2=3, \dim \fl^3=0$, or\vspace{-1ex} 
\item[(iv)] $\dim \fl^2=3, \dim \fl^3=2, \dim \fl^4=0$. 
\end{itemize} 
If (i) holds, then $\fl\cong \RR^5$. If (ii) holds, then $\fl\cong \fg_{5,2}$ 
by 
Proposition \ref{P2}. The conditions in (iii) cannot be satisfied for a 
5-dimensional Lie 
algebra $\fl$, since in this case $\fl'\subset \fz(\fl)$, thus $\dim \fz \ge 3$ 
and therefore $\dim \fl'\le 1$, which contradicts $\dim \fl'=3$. 

Now assume that (iv) holds. Choose linear independent vectors $X_1,X_2$ in 
$\fl\setminus 
\fl'$. Then $X_3:=[X_1,X_2]\notin \fl^3$. For $X_4:=[X_1,X_3]$ and 
$X_5:=[X_2,X_3]$ we now have 
$\fl^3=[\fl',\fl]=[X_3+\fl^3,\fl]=[X_3,\fl]=\Span \{X_4,X_5\}$. Hence, 
$X_1,\dots,X_5$ is a 
basis of $\fl$. Since $X_4,X_5\in\fl^3$ are central we obtain 
$$\begin{array}{llrrrl} 
0\ =\ d\gamma(X_1,X_2,X_3,X_4)&=&\gamma(X_1,[X_2,X_3],X_4)& =& 
\gamma(X_1,X_5,X_4)&\\[0.5ex] 
0\ =\ d\gamma(X_1,X_2,X_3,X_5)&=&\gamma([X_1,X_3],X_2,X_5)& =& 
\gamma(X_4,X_2,X_5)&\\[0.5ex] 
0\ =\ d\gamma(X_1,X_2,X_4,X_5)&=&-\gamma([X_1,X_2],X_4,X_5)& =& 
-\gamma(X_3,X_4,X_5)&\hspace{-0.5em}, 
\end{array}$$ 
thus $\gamma(X_4,X_5,\cdot)=0$, which contradicts Condition $(A_2)$. 
\qed 
\begin{pr} \label{P5} 
\begin{enumerate} 
\item 
If $\fl=\RR^5$ and $\fa=0$, then $\cHQ_0/G$ consists of one element. This 
element is 
represented by $[0,\gamma_0]\in \cHQ_0$, where 
$\gamma_0=(\sigma^1\wedge\sigma^2 +\sigma^3\wedge\sigma^4)\wedge \sigma^5$ 
for a fixed basis $\sigma^1,\dots,\sigma^5$ of $\fl^*$. 
\item 
If $\fl=\fg_{5,2}=\{[X_1,X_2]=Y,\ [X_1,X_3]=Z \}$ and $\fa=0$, then 
$\cHQ_0/G$ 
consists of 
two elements. These elements are represented by $[0,\gamma_1], [0,\gamma_2]\in 
\cHQ_0$, 
where 
$\gamma_1=\sigma^1\wedge\sigma^Y\wedge \sigma^Z$ and 
$\gamma_2=\sigma^1\wedge\sigma^Y\wedge 
\sigma^Z +\sigma^2\wedge\sigma^3\wedge \sigma^Z$ 
for the basis $\sigma^1,\sigma^2,\sigma^3,\sigma^Y, \sigma^Z$ of $\fl^*$ which 
is 
dual to 
$X_1,X_2,X_3,Y,Z$. 
\end{enumerate} 
\end{pr} 
\proof The statement for $\fl=\RR^5$ is easy to prove. Take $\fl=\fg_{5,2}$. For 
$c\in\RR\setminus 0$, $A\in \GL(2,\RR)$, $y =(y_1,y_2)$, $y_1,y_2\in\RR^2$, 
$x\in 
\gl(2,\RR)$ we 
define a linear map $S(c,A,x,y):\fg_{5,2}\rightarrow\fg_{5,2}$ by 
$$S(c,A,x,y)= \mbox{{\small $
\left(\begin{array}{ccc} 
c&0&0\\ 
y_1& A&0\\ 
y_2&x&cA 
\end{array}\right)$}}$$ 
with respect to the basis $X_1,X_2,X_3,Y,Z$ of $\fg_{5,2}$. Using that 
$\fq:=\Span\{X_2,X_3,Y,Z\}$ is the unique 4-dimensional abelian ideal of 
$\fg_{5,2}$ 
and 
that $\Span\{Y,Z\}$ is the centre of $\fg_{5,2}$ it is not hard to show that the 
automorphism group of $\fg_{5,2}$ equals 
$$\Aut (\fg_{5,2})= \{S(c,A,x,y)\mid c\in\RR\setminus 0, A\in GL(2,\RR), 
y_1,y_2\in\RR^2, 
x\in \gl(2,\RR) 
\}\,. 
$$   
Obviously, we have $\cH^2_Q(\fl,0)=H^3(\fg_{5,2})$. Using the Hochschild-Serre 
spectral 
sequence we see that $H^3(\fg_{5,2})$ is determined by the exact sequence 
$$ 
0\longrightarrow H^1(\RR\cdot X_1,H^2(\fq))\longrightarrow 
H^3(\fg_{5,2})\longrightarrow 
H^0(\RR\cdot X_1,H^3(\fq))\longrightarrow 0\,. 
$$  
We have 
$$ H^1(\RR\cdot X_1,H^2(\fq))=C^1(\RR\cdot X_1,C^2(\fq))/B^1(\RR\cdot 
X_1,C^2(\fq)),$$ 
where 
$$B^1(\RR\cdot X_1,C^2(\fq))=\left\{\sigma\in C^1(\RR\cdot X_1,C^2(\fq))\ \Big| 
\begin{array}{l} \sigma(X_1)(X_2,Z)+\sigma(X_1)(X_3,Y)=0\\ 
\sigma(X_1)(X_2,Y)=\sigma(X_1)(X_3,Z)=0\\\sigma(X_1)(Y,Z)=0 
\end{array}\right\}\,,$$ 
and 
\begin{equation} \label{EH}
H^0(\RR\cdot X_1,H^3(\fq))=C^3(\fq)^{ X_1} 
=\{ \sigma\in C^3(\fq) \mid \sigma(X_2,Y,Z)=\sigma(X_3,Y,Z)=0\}\,. 
\end{equation} 

Observe that $\cH^2_Q(\fl,0)_0=\cH^2_Q(\fl,0)_\sharp$ since $\fl$ is not the 
direct sum of two non-trivial Lie algebras. In particular, Condition $(A_1)$ and 
Equation (\ref{EH}) imply
$$\cH^2_Q(\fl,0)_0=\{[\gamma]\in H^3(\fg_{5,2})\mid 
\gamma(X_1,Y,Z)\not=0\}.$$ Using the description of $ H^1(\RR\cdot 
X_1,H^2(\fq))$ 
given above 
we see that 
$$\{S(c,\Id,x,0)\mid c\in\RR\setminus 0,x\in \gl(2,\RR)\}\subset \Aut 
(\fg_{5,2})$$ 
acts transitively on $\{[\sigma]\in H^1(\RR\cdot X_1,H^2(\fq)) \mid 
\sigma(X_1)(Y,Z)\not=0\}$. Furthermore, using the description of $ H^0(\RR\cdot 
X_1,H^3(\fq))$ we see that the action of 
$$\{S(1,A,0,0)\mid \det A=1\}\subset\Aut (\fg_{5,2})$$ 
on $ H^0(\RR\cdot X_1,H^3(\fq))$ 
has two orbits represented by $\sigma_1=0$ and 
$\sigma_2=\sigma^2\wedge\sigma^3\wedge 
\sigma^Z$. Moreover, this group leaves $\sigma^1\wedge \sigma^Y \wedge \sigma^Z$ 
invariant. 

It is easy to check that the orbits of $[0,\gamma_1]$ and $[0,\gamma_2]$ are 
different.
\qed

\begin{pr}\label{g41} Take $\fl=\fg_{4,1}=\{[X_1,Z]=Y,\ [X_1,X_2]=Z\}$ and let 
$\fa$ 
be a 
trivial 
$\fl$-module. If $\fa=0$ or $\dim 
\fa \ge3$, then $\cHQ_0=\emptyset$. 

If $\dim \fa=1$, then $\cHQ_0/G$ consists of 
two elements. They are represented by 
$[\alpha,\gamma]=[\sigma^{1}\wedge\sigma^{Y}\otimes A,0]$ and $[\alpha,\gamma]=[\sigma^{1}\wedge\sigma^{Y}\otimes A, 
\sigma^{2}\wedge\sigma^{Y}\wedge\sigma^{Z}].$
where $A$ is a fixed unit 
vector in 
$\fa$. 

If $\fa\in\{\RR^2,\RR^{2,0}\}$, then $\cHQ_0/G$ consists of 
one element, which is represented by  
\begin{equation}\label{qu1}
[\alpha,\gamma]=[\sigma^{1}\wedge\sigma^{Y}\otimes A_{1}+ 
 \sigma^{2}\wedge\sigma^{Z}\otimes A_{2}, 
 0]\end{equation} 
where $A_{1},A_{2}$ is a fixed orthonormal basis in $\fa$. 
If $\fa=\RR^{1,1}$, then $\cHQ_0/G$ also consists of 
two elements, which are represented also by  
(\ref{qu1}) but where now either $A_1,A_2$ or $A_2,A_1$ is an 
orthonormal basis.
\end{pr} 
\proof 
Let us first determine the automorphism group of $\fl$. For $a,b,c\in\RR$ and 
$x=(x_1,\dots,x_4)\in\RR^4$ we define a linear map $S(a,b,c,x):\fl\rightarrow 
\fl$ by 
$$ 
S(a,b,c,x)=\mbox{{\small
$\left(\begin{array}{cccc} 
a&0&0&0\\x_1&b&0&0\\x_2&c&ab&0\\ 
x_3&x_4&ac&a^2b 
\end{array} 
\right)$}} 
$$ 
with respect to the basis $X_1,X_2,Z,Y$ of $\fl$. Using that the unique 
3-dimensional 
abelian ideal $\fq=\Span\{X_2,Y,Z\}$ of $\fl$, $\fl'=\Span\{Y,Z\}$ and 
$\fl^3=\RR\cdot Y$ 
are invariant under each automorphism of $\fl$ it is not hard to check that the 
automorphism group of $\fg_{4,1}$ equals 
$$\Aut (\fg_{4,1})= \{S(a,b,c,x)\mid a,b\in\RR\setminus 0, c\in\RR, x\in 
\RR^4\}\,.$$ 

The cohomology group $H^2(\fg_{4,1},\fa)$ is determined by the exact sequence 
$$ 
0\longrightarrow H^1(\RR\cdot X_1,H^1(\fq,\fa))\longrightarrow 
H^2(\fg_{4,1},\fa)\longrightarrow H^0(\RR\cdot X_1,H^2(\fq,\fa))\longrightarrow 
0\,. 
$$  
We have 
$$ H^1(\RR\cdot X_1,H^1(\fq,\fa))=C^1(\RR\cdot X_1,C^1(\fq,\fa))/B^1(\RR\cdot 
X_1,C^1(\fq,\fa)),$$ 
where 
$$B^1(\RR\cdot X_1,C^1(\fq,\fa))=\{\sigma\in C^1(\RR\cdot X_1,C^1(\fq,\fa))\mid 
\sigma(X_1)(Y)=0\}\,,$$ 
and 
\begin{eqnarray*} 
H^0(\RR\cdot X_1,H^2(\fq,\fa))&=&C^2(\fq,\fa)^{ X_1} 
=\{ \sigma\in C^2(\fq,\fa) \mid \sigma(X_2,Y)=\sigma(Y,Z)=0\}\,. 
\end{eqnarray*} 
In particular, $(A_2)$ implies $\alpha(Y,X_1)\not= 0$. If $\alpha(Y,X_1)\not= 
0$, 
then also $(A_0)$ and $(A_1)$ are satisfied. Since $d\gamma=0$ for all 
$\gamma\in C^3(\fl)$ the equation $2d\gamma=\langle 
\alpha\wedge \alpha\rangle $ holds if and only if 
$\alpha(Y,X_1)\perp\alpha(X_2,Z)$. 
Hence we obtain 
$$\cHQ_0=\left \{[\alpha,\gamma]\in \cHQ\ \Big| 
\begin{array}{l} 
\alpha= (\sigma^1\wedge \sigma^Y)\otimes A_1+ (\sigma^2\wedge \sigma^Z)\otimes 
A_2, 
\\ 
\fa=\Span\{A_1,A_2\},\ A_1\perp A_2,\ A_1\not=0 
\end{array} \right\}\,.$$ 
 
Because of $B^3(\fl)=\{\gamma\in C^3(\fl)\mid 
\gamma(X_1,Y,Z)=\gamma(X_2,Y,Z)=0\}$ we 
may assume that 
$$\gamma= r\sigma^2\wedge \sigma^Y \wedge \sigma^Z+s\sigma^1\wedge \sigma^Y 
\wedge 
\sigma^Z\,. $$ 

Next we discuss the action of $G$ on ${\cal H}_Q^2(\fl,\fa)_0$. Let $b\in\RR$ be such that $\langle A_1,A_1\rangle=\kappa b^{-2}$, $\kappa=\pm1$, and put $S=S(1,b,0,x)$, $x=(0,0,0,-\kappa s b^3)$. Moreover, set $\tau(Z):=-\kappa s b^3 A_1$, $\tau(X_1)=\tau(X_2)=\tau(Y)=0$. Then $(S^*\gamma +\langle (S^*\alpha+\textstyle{\frac 12} d\tau)\wedge \tau\rangle)(X_1,Y,Z)=0$. Hence there exists $\sigma\in C^3(\fl)$ such that $(S^*(\alpha,\gamma))\cdot (\tau,\sigma)= ((\sigma^1\wedge \sigma^Y)\otimes b A_1+(\sigma^2\wedge \sigma^Z)\otimes A_2,\hat r \sigma^2\wedge \sigma^Y\wedge \sigma^Z)$, $\hat r \in \RR$. This shows that we may assume $s=0$ and $\langle A_1,A_1 
\rangle =\pm 1$ without changing the $G$-orbit.

Now let us assume $\dim \fa=1$. Take 
$[\alpha,\gamma]\in 
\cHQ_0$, $\alpha= (\sigma^1\wedge \sigma^Y)\otimes A_1$, $\gamma= r\sigma^2\wedge \sigma^Y\wedge \sigma^Z$.
Suppose $r\not=0$. Then $r=\kappa r_+$, $\kappa=\pm 1$, $r_+>0$. Choose $a,b\in\RR$ such that $b^2=1/r_+$ and $a^3=\kappa/b$  holds. Then $S:=S(a,b,0,0)$ satisfies $(S,\kappa\cdot\Id_\fa)^*[\alpha,\gamma]=[\alpha,\sigma^2\wedge\sigma^Y\wedge\sigma^Z]$. Furthermore, $[\alpha,\sigma^2\wedge\sigma^Y\wedge\sigma^Z]$ and $[\alpha,0]$ are not in the same $G$-orbit since 
$$\langle S^*\alpha\wedge\tau\rangle(X_2, Y, Z)=\langle d\tau\wedge \tau\rangle (X_2, Y, Z)=d\sigma(X_2, Y, Z)=0$$
for all $(S,U)\in G$ and $(\tau,\sigma)\in {\cal C}^1_Q(\fl,\fa)$.

Consider now the case $\dim \fa=2$. Take 
$\alpha= (\sigma^1\wedge \sigma^Y)\otimes A_1+ (\sigma^2\wedge 
\sigma^Z)\otimes 
A_2$ and  $\gamma= r\sigma^2\wedge \sigma^Y \wedge \sigma^Z $. If  $\langle A_2,A_2\rangle =\pm c^{-2}$, then choose 
$a,b$ such that $a^3b=1$, $ab^2=c$. Applying
$S(a,b,0,0)$ to 
$[\alpha,\gamma]$ we see that we may assume $\langle A_2,A_2\rangle =\pm 1$. 
Thus $\langle A_i,A_i\rangle =\kappa_i$, $i=1,2$. Put $S:=S(1,1,0,x)$ for $x=(r\kappa_2,0,0,-r^2\kappa_1\kappa_2/2)$ and define $\tau\in C^1(\fl,\fa)$ by $\tau(X_1)=\tau(X_2)=0$, $\tau(Y)=r\kappa_2 A_2$, $\tau(Z)=-r^2\kappa_1\kappa_2/2$. Then $(S^*(\alpha,\gamma))\cdot(\tau,\sigma)=(\alpha,0)$ for suitable $\sigma\in C^3(\fl)$.
\qed 

Now let $\fl$ be one of the Lie algebras $\fh(1)\oplus\RR= 
\{[X_1,X_2]=X_3\}\oplus \RR\cdot X_4$ or $\RR^4=\Span\{X_1,\dots,X_4\}$. 
Let $\sigma^1,\dots,\sigma^4$ be a basis of $\fl^*$ which is 
dual to 
$X_1,\dots,X_4$. Let $A_1,A_2,\ldots$ be a basis of a vector space $\fa$. We 
define the 
following 
2-forms 
$$\begin{array}{lll} 
\alpha_1 &=& (\sigma^1\wedge\sigma^3 +\sigma^2\wedge\sigma^4)\otimes A_1+ 
(\sigma^2\wedge\sigma^3 +\sigma^1\wedge\sigma^4)\otimes A_2\\[0.3ex] 
\alpha_2 &=& (\sigma^1\wedge\sigma^3 -\sigma^2\wedge\sigma^4)\otimes A_1+ 
(\sigma^2\wedge\sigma^3 +\sigma^1\wedge\sigma^4)\otimes A_2\\[0.3ex] 
\alpha_3 &=& (\sigma^1\wedge\sigma^3)\otimes A_1+ 
(\sigma^2\wedge\sigma^3 +\sigma^1\wedge\sigma^4)\otimes A_2\\[0.3ex] 
\alpha_4 &=& (\sigma^1\wedge\sigma^3)\otimes A_1+ 
(\sigma^2\wedge\sigma^3)\otimes A_2\\[0.3ex] 
\alpha_5 &=& (\sigma^1\wedge\sigma^3)\otimes A_1+ 
(\sigma^1\wedge\sigma^4)\otimes A_2,
\quad \alpha_5'= (\sigma^1\wedge\sigma^4)\otimes A_1+ 
(\sigma^1\wedge\sigma^3)\otimes A_2\\[0.3ex] 
\alpha_6 &=& (\sigma^1\wedge\sigma^3)\otimes A_1+ 
(\sigma^2\wedge\sigma^4)\otimes A_2,
\quad \alpha_6'= (\sigma^2\wedge\sigma^4)\otimes A_1 
+(\sigma^1\wedge\sigma^3)\otimes A_2\\[0.3ex] 
\alpha_7 &=& (\sigma^1\wedge\sigma^3)\otimes A_1 \,.
\end{array} 
$$ 
Moreover, we define the 3-form $\gamma_0$ on $\fl$ by 
$\gamma_0=\sigma^2\wedge\sigma^3\wedge\sigma^4$.

\begin{pr} \label{Lh1} Take $\fl=\fh(1)\oplus\RR= 
\{[X_1,X_2]=X_3\}\oplus \RR\cdot X_4$. Let $\fa$ be a trivial 
orthogonal 
$\fl$-module. 
If $\fa=\RR^2$ or $\fa=\RR^{2,0}$, then the elements in $\cHQ_0/G$ are 
represented by 
$[\alpha_1,0]$, $[\alpha_5,0]$, $[\alpha_5,\gamma 
_0]$, $[\alpha_6,0]$, $[\alpha_6,\gamma 
_0]$, where $A_1,A_2$ is a fixed orthonormal basis of~$\fa$. 

If $\fa=\RR^{1,1}$, then $\cHQ_0/G$ has eleven elements, three of them are 
represented by 
$[\alpha_1,0]$, $[\alpha_2,0]$, $[\alpha_3,0]$, where $A_1,A_2$ is a fixed Witt 
basis 
of 
$\fa$, eight further elements are represented by $[\alpha_5,0]$, 
$[\alpha_5,\gamma_0]$, 
$[\alpha_6,0]$, $[\alpha_6,\gamma_0]$, $[\alpha_5',0]$, $[\alpha_5',\gamma_0]$, 
$[\alpha_6',0]$, $[\alpha_6',\gamma_0]$, where $A_1,A_2$ is a fixed orthonormal 
basis of $\fa$. 

If $\fa=\RR^1$ or $\fa=\RR^{1,0}$, then there is only one element in $\cHQ_0/G$. 
It is represented by $[\alpha_7,\gamma_0]$, where $A_1$ is a fixed unit vector 
in 
$\fa$.

If $\fa=0$, then $\cHQ_0=\emptyset$.
\end{pr} 

\proof 
For $A,X\in \gl(2,\RR)$ and $u=(u_1,u_2,u_3)\in\RR^3$ we define 
$$S(A,X,u)=\mbox{{\small $\left( \begin{array}{cc} A &0\\X & U \end{array} 
\right)$}}\in\gl(4,\RR), 
\mbox{ 
where } U= 
\mbox{{\small $\left( \begin{array}{cc} u_1 &u_2\\0 & u_3 \end{array} \right)$}} 
\in\gl(2,\RR).$$ 
Then the automorphism group of $\fl=\fh(1)\oplus\RR$ equals 
$$\Aut(\fl)=\left\{ S(A,X,u)\mid u\in\RR^3,\ X\in\gl(2,\RR),\ A\in\GL(2,\RR),\ 
\det 
A=u_1,\ 
u_3\not=0 \right\}, $$ 
where we consider all automorphisms with respect to the basis $X_1,\dots,X_4$ of 
$\fl$. 

By direct computations or using the K\"unneth formula and the explicit 
description of $H^2(\fh(1),\fa)$ in \cite{KO2}, we see that 
\begin{eqnarray*} 
Z_\fl:=\{\alpha\in C^2(\fl,\fa)\mid \alpha(X_1,X_2)=\alpha(X_3,X_4)=0\} 
&\longrightarrow& 
H^2(\fl,\fa)\\ 
\alpha &\longmapsto& [\alpha] 
\end{eqnarray*} 
is a bijection. 

Now take $[\alpha,\gamma]\in\cHQ_0$, $\alpha\in Z_\fl$. Since, obviously, 
$d\gamma=0$ we have 
${\langle \alpha\wedge\alpha \rangle =0}$. 
Condition $(A_1)$ gives $\alpha(X_3,\fl)\not=0$ and Condition $(B_1)$ says that 
$\alpha(X_3,\fl)$ is non-degenerate. By indecomposability we have 
$\alpha(\fl,\fl)=\fa$. 
Hence $\alpha$ is an element of the $G$-invariant subset 
$C\subset Z_\fl$ defined by 
$$C:=\{\alpha \in Z_\fl \mid \langle \alpha\wedge\alpha \rangle =0,\ 
\alpha(\fl,\fl)=\fa,\ 
0\not= \alpha(X_3,\fl)\subset \fa \mbox{ is non-degenerate} \}.$$ 

A cocycle $\alpha\in Z_\fl$ satisfies 
$\langle\alpha\wedge\alpha\rangle=0$ if and only if 
\begin{equation}\label{pair} 
\langle \alpha(X_1,X_3),\alpha(X_2,X_4)\rangle =\langle 
\alpha(X_2,X_3),\alpha(X_1,X_4)\rangle. 
\end{equation} 

Let us determine the $G$-orbits in $C$ in the case that $\dim \fa\le 2$. Take 
$\alpha\in C$. 
In particular we have $\dim \alpha(X_3,\fl)=1$ or $\dim \alpha(X_3,\fl)=2$.

Let us first consider the case $\dim \alpha(X_3,\fl)=1$. Replacing $\alpha$ by 
an 
element in 
the $G$-orbit of $\alpha$ we may assume $\alpha(X_1,X_3)=A_1$ and 
$\alpha(X_2,X_3)=0$, where 
$\langle A_1,A_1 \rangle =\pm1$. From (\ref{pair}) we obtain $\langle 
\alpha(X_1,X_3),\alpha(X_2,X_4)\rangle =0$. Hence either $\alpha(X_2,X_4)=0$ or 
$\alpha(X_2,X_4)=A_2\not=0$ and $A_1\perp A_2$. In the latter case $A_2$ cannot 
be 
isotropic 
since $\fa$ is at most two-dimensional. Hence, replacing $\alpha$ by an element 
in the same 
$G$-orbit we may assume that $\langle A_2,A_2\rangle=\pm1$  and that 
$\alpha(X_1,X_4)=0$, hence $\alpha$ is in the same $G$-orbit as $\alpha_6$ or as 
$\alpha_6'$ for an orthonormal basis $A_1,A_2$. In 
the first case, where $\alpha(X_2,X_4)=0$ we may assume that $\alpha(X_1,X_4)=0$ 
or 
$\alpha(X_1,X_4)=A_2\not=0$ and $A_1\perp A_2$, $\langle A_2,A_2 \rangle =\pm1$, 
hence $\alpha$ is in the same orbit as $\alpha_5$ or $\alpha_5'$ for an 
orthonormal basis $A_1,A_2$ or as $\alpha_7$ for a unit vector $A_1$. Obviously, 
the orbit 
of $\alpha_7$ contains neither $\alpha_5$, $\alpha_5'$, $\alpha_6$, nor 
$\alpha_6'$. Also the orbits of $\alpha_5$ and $\alpha_6$ are different. Indeed, 
$\alpha_5(X_2,\fl)=0$ and $\alpha_6(L,\fl)\not=0$ for all $L\in\fl, L\not=0$. 
Analogously, the orbits of $\alpha_5'$ and $\alpha_6'$ are different. Moreover, 
$\alpha_i$ and $\alpha_i'$, $i=5,6$, are not on the same orbit, since $X_3$ 
plays a distinguished role in $\fl$.

Now we consider the case $\dim \alpha(X_3,\fl)=2$. Then $\alpha(X_1,X_3)=:A_1$ 
and 
$\alpha(X_2,X_3)=:A_2$ are linearly independent. First we show that we may 
assume 
that $A_1, 
A_2$ is an orthonormal basis of $\fa$ if $\fa=\RR^2$ or $\fa=\RR^{2,0}$ and that 
$A_1, A_2$ 
is a Witt basis if $\fa=\RR^{1,1}$ (replacing $\alpha$ by an element in the same 
$G$-orbit). 
Clearly, we can choose $X_1',X_2'\in \Span\{X_1,X_2\}$ such that  
$\alpha(X_1',X_3),\, 
\alpha(X_2',X_3)$ is an orthonormal basis or a Witt basis, respectively. We have 
$[X_1',X_2']=rX_3$, $r\in\RR$, $r\not=0$. We choose $s\in\RR$ such that $s^3=r$ 
and 
define 
$\bar X_1=(1/s) X_1'$, $\bar X_2=(1/s) X_2'$, and $\bar X_3=sX_3$. Take $\bar 
X_1, 
\bar X_2, 
\bar X_3, X_4$ as a new basis for $\fl$. Hence, we may assume that $A_1, A_2$ 
are as 
claimed. By (\ref{pair}) we have 
\begin{equation}\label{pair1} 
\langle A_1,\alpha(X_2,X_4)\rangle =\langle A_2,\alpha(X_1,X_4)\rangle. 
\end{equation} 
Replacing $X_4$ by a suitable linear combination of $X_4$ and $X_3$ we may 
assume 
that  
$\alpha(X_2,X_4)$ is a multiple of $A_1$.

Assume that $\fa=\RR^2$ or $\fa=\RR^{2,0}$. If $\alpha(X_2,X_4)=0$, then 
(\ref{pair1}) 
implies that $\alpha(X_1,X_4)$ is a multiple of $A_1$. Hence we may assume that 
either 
$\alpha(X_1,X_4)=0$ or that $\alpha(X_1,X_4)=A_1$. 
Consequently, 
$\alpha$ is in the same orbit as $\alpha_4$ or as the 2-form 
$$\alpha_1':=(\sigma^1\wedge\sigma^3 
+\sigma^1\wedge \sigma^4)\otimes A_1 +(\sigma^2\wedge\sigma^3)\otimes A_2,$$ 
which is in the same orbit as $\alpha_1$. Indeed, we have $U^{-1}\circ 
S^*\alpha_1=\alpha_1'$ for $U$, $S=S(A,0,u)$ with
$$U=A=\mbox{$\frac 1{\sqrt 2}${\small 
$\left(\begin{array}{cc}1&-1\\1&1\end{array}\right)$,}} \quad 
u=(1,\,1/2,\,1/2),$$ 
where we take $U$ with respect to the basis $A_1,A_2$ of $\fa$.

If $\alpha(X_2,X_4)=rA_1$, $r\not=0$, then rescaling $X_4$ we may assume that 
$\alpha(X_2,X_4)=A_1$. Now (\ref{pair1}) yields $\alpha(X_1,X_4)= A_2 +sA_1$. We 
will 
show 
that we may assume $s=0$. For $s\in \RR$ we choose $t\in \RR$ such that $s=2 
\tan 2t$ 
and we 
define 
$$a=\sin t,\ b=\cos t, \ u_2= \sin 2t,\ u_3=-\cos 2t$$ 
and $$ A=\mbox{{\small $\left( \begin{array}{cc} a &-b\\b & a \end{array} 
\right)$}}, \ 
X=0\in\gl(2,\RR), \ 
u=(1,u_2,u_3)\in\RR^3.$$ 
For this choice of $A,X$, and $u$ we consider $S=S(A,X,u)\in\Aut(\fl)$ and we 
define 
$X_i'=SX_i$, $i=1,\dots,4$. Then we have 
\begin{eqnarray*} 
\alpha(X_1', X_3')&=& aA_1 +bA_2 \,=:\, A_1'\\ 
\alpha(X_2', X_3')&=& -bA_1 +aA_2 \,=:\, A_2'\\ 
\alpha(X_1', X_4')&=& (u_2a+au_3s +bu_3)A_1 +(u_2b+au_3)A_2 \,=\, A_2'\\ 
\alpha(X_2', X_4')&=& (-u_2b-sbu_3+au_3)A_1 +(u_2a-bu_3)A_2 \,=\, A_1', 
\end{eqnarray*} 
where $A_1',A_2'$ is again an orthonormal basis. Hence, $\alpha$ is in the same 
orbit 
as 
$\alpha_1$. The 2-forms $\alpha_1$ and $\alpha_4$ are on different orbits, since 
$\alpha_4(X_4,\fl)=0$ and $\alpha_1(L,\fl)\not=0$ for all $L\in\fl, L\not=0$. 

Take now $\fa=\RR^{1,1}$. Recall that we may assume $\alpha(X_1,X_3)=A_1$ and 
$\alpha(X_2,X_3)=A_2$ such that $A_1, A_2$ is a Witt basis and that 
$\alpha(X_2,X_4)=rA_1$, 
$r\in\{0,1\}$. From (\ref{pair1}) we get $\alpha(X_1,X_4)=s A_2$ for a real 
number 
$s$. If 
$r=s=0$, then $\alpha$ is in the same orbit as $\alpha_4$. If $r=0$, $s=1$ or 
$r=1$, 
$s=0$, 
then $\alpha$ is in the same orbit as $\alpha_3$. If $r=1$, $s\not=0$, then we 
put 
$x=|s|^{-1/4}$ and $v=(\sgn s)\cdot |s|^{-1/2}$. We define 
$S=\diag(x,x^{-1},1,v)\in 
\Aut(\fl)$ and $U=\diag(x^{-1},x)$. Then $(S,U)^*\alpha$ equals $\alpha_1$ or 
$\alpha_2$. 
The 2-forms $\alpha_1,\dots,\alpha_4$ are on different orbits, since the 
elements 
of 
its 
orbits differ in the properties of their projections to the isotropic lines in 
$\fa=\RR^{1,1}$. 

We can summarize this as follows. 
If $\fa=\RR^2$ or $\fa=\RR^{2,0}$, then there are four $G$-orbits in $C$ 
represented 
by 
$\alpha_1, \alpha_4, \alpha_5, \alpha_6$, where $A_1,A_2$ is a fixed 
orthonormal basis of $\fa$. 
If $\fa=\RR^{1,1}$, then there are eight $G$-orbits in $C$, four of them are 
represented by 
$\alpha_1, \alpha_2, \alpha_3, \alpha_4$, where $A_1,A_2$ is a fixed Witt 
basis 
of $\fa$, four further orbits are represented by $\alpha_5,\alpha_5', \alpha_6, 
\alpha_6'$, where now
$A_1,A_2$ is 
a fixed orthonormal basis of $\fa$. 
If $\fa=\RR^1$ or $\fa=\RR^{1,0}$, then $\alpha_7\in C$ and $G$ acts 
transitively on 
$C$.

Since $Z^1(\fl,\fa)=\{\tau\in C^1(\fl,\fa) \mid d\tau=0\}=\{\tau\in C^1(\fl,\fa) 
\mid 
\tau(X_3)=0\}$ we have $\langle \alpha_i \wedge Z^1(\fl,\fa) \rangle =C^3(\fl)$ 
for 
$i=1,\dots,3$. For $\alpha_4$ we have $B^3(\fl)+\langle \alpha_4 \wedge 
Z^1(\fl,\fa) 
\rangle =C^3(\fl)$, where $B^3(\fl)=\{d\sigma\mid \sigma \in 
C^2(\fl)\}=\RR\cdot\sigma^1\wedge \sigma^2\wedge\sigma^4$. Hence 
$[\alpha_i,\gamma]=[\alpha_i,0]\in\cHQ$ for $i=1,\dots,4$ and all 
$[\alpha_i,\gamma]\in\cHQ$. 
Note that $[\alpha_4,0]$ is decomposable. 

If $\alpha\in\{\alpha_5,\alpha_5',\alpha_6,\alpha_6',\alpha_7\}$, then 
$\gamma_0$ spans a complement of 
$B^3(\fl)+\langle \alpha \wedge Z^1(\fl,\fa) \rangle$ in $C^3(\fl)$. Hence, 
for all these $\alpha$ and for all $\gamma\in C^3(\fl)$ there exists a real 
number $c$ such that 
$[\alpha,\gamma]=[\alpha,c\gamma_0]\in\cHQ$. Let us first determine the 
$G$-orbit of 
$[\alpha_i,c\gamma_0]$ for $i=6,7$. If $c\not=0$, then 
$S:=\diag(c,1/c^2,1/c,c^2)\in\Aut(\fl)$ and we have 
$(S,\Id)^*[\alpha_i,c\gamma_0]=[\alpha_i,\gamma_0]$. Obviously the orbits of 
$[\alpha_i,\gamma_0]$ and $[\alpha_i,0]$ are different. Now  
consider the $G$-orbit 
of $[\alpha_5,c\gamma_0]$. 
If $c\not=0$ we put $s=|c|^{1/4}$. Then $S_5:=\diag(s, 
1/s^2,1/s,(\sgn 
c)\cdot(1/s))\in\Aut(\fl)$, $U:=\diag(1,\sgn c)\in O(\fa)$ and we get 
$(S_5,U)^*[\alpha_5,c\gamma_0]=[\alpha_5,\gamma_0]$. The orbits of 
$[\alpha_5,\gamma_0]$ and $[\alpha_5,0]$ are different. Analogously, one 
determines the orbits of $[\alpha_i',\gamma]$ for $i=5,6$.

It remains to check admissibility and indecomposability. All cohomology classes 
$[\alpha,\gamma]\in\cHQ$ with $\alpha\in C$ satisfy $(B_0)$, $(A_1)$, and 
$(B_1)$. Moreover, it is not hard to see that all cohomology classes listed in 
the proposition satisfy also $(A_0)$ and are indecomposable. 
\qed

\begin{pr} \label{Pl} Let $\fl$ be the abelian Lie algebra 
$\RR^4=\Span\{X_1,\dots,X_4\}$ and let $\fa$ be a trivial orthogonal 
$\fl$-module. 

If $\fa=\RR^2$ or $\fa=\RR^{2,0}$, then the elements in $\cHQ_0/G$ are 
represented by 
$[\alpha_1,0]$ and $[\alpha_4,\sigma^1\wedge\sigma^2\wedge\sigma^4]$, where 
$A_1,A_2$ 
is a 
fixed orthonormal basis of $\fa$. 

If $\fa=\RR^{1,1}$, then the elements in $\cHQ_0/G$ are represented by 
$[\alpha_1,0]$, 
$[\alpha_2,0]$, $[\alpha_3,0]$, 
$[\alpha_4,\sigma^1\wedge\sigma^2\wedge\sigma^4]$
where $A_1,A_2$ is a fixed Witt basis of $\fa$. 

If $\fa=\RR^1$ or $\fa=\RR^{1,0}$, then $\cHQ_0/G$ contains exactly one element. 
This 
is represented by $[\alpha_7,\gamma_0]$, where $A_1$ is a fixed unit vector in 
$\fa$. 

If $\fa=0$, then $\cHQ_0=\emptyset$.
\end{pr} 
\proof First notice that $[0,\gamma]\in\cHQ$ is in the same $G$-orbit as either 
$[0,0]$ or $[0,\gamma_0]$. Since both of these cohomology classes are 
decomposable $[\alpha,\gamma]\in\cHQ_0$ implies $\alpha\not=0$. Hence, under the 
assumptions of the proposition we have $\dim\alpha(\fl,\fl)=1$ or 
$\dim\alpha(\fl,\fl)=2$.

If  $\dim \alpha(\fl,\fl)=2$, then we may assume that $\alpha(X_1,X_3)$ and  
$\alpha(X_2,X_3)$ are linearly independent and that $\alpha(X_1,X_2)=0$.
This can easily be verified using the same idea as in the proof of (\ref{E10n}).
Hence, Equation (\ref{pair}) also holds in this case
and we can argue as in the proof of Proposition~\ref{Lh1} to show that 
$\alpha$ is 
in the same $G$-orbit as $\alpha_1$ or $\alpha_4$ if $\fa=\RR^2$ or 
$\fa=\RR^{2,0}$ 
and in the same orbit as one of the 2-forms $\alpha_1,\dots,\alpha_4$ if 
$\fa=\RR^{1,1}$ and that all these orbits are different. 

If  $\dim \alpha(\fl,\fl)=1$, then by classification of ordinary 2-forms there 
exists a map $S\in\Aut(\fl)=\GL(4,\RR)$ such that  $S^*\alpha=\alpha_7$ or 
$S^*\alpha=\alpha':=(\sigma^1\wedge\sigma^3 +\sigma^2\wedge\sigma^4)\otimes 
A_1$. Since  $\langle \alpha'\wedge \alpha'\rangle\not=0$ we can exclude the 
latter case.
 
Again we have $\langle \alpha_i\wedge Z^1(\fl,\fa)\rangle =C^3(\fl)$ for 
$i=1,2,3$. 
Furthermore, $\RR\cdot \sigma^1\wedge \sigma^2\wedge \sigma^4$ is a complement 
of  $\langle \alpha_4\wedge Z^1(\fl,\fa)\rangle $ in $C^3(\fl)$ and $\Span\{ 
\sigma^1\wedge \sigma^2\wedge \sigma^4,\gamma_0\}$ is a complement of  
$\langle \alpha_7\wedge Z^1(\fl,\fa)\rangle $ in $C^3(\fl)$. 
Note that $[\alpha_4,0]$ and $[\alpha_7,0]$ are decomposable. Hence, if 
$[\alpha_4,\gamma]\in\cHQ_0$, then $[\alpha_4,\gamma]$ is on the same $G$-orbit  
as $[\alpha_4,\sigma^1\wedge\sigma^2\wedge\sigma^4]$. Moreover, if 
$[\alpha_7,\gamma]\in\cHQ_0$, then $[\alpha_7,\gamma]= [\alpha_7,\gamma']$, 
where $\gamma'=\sigma^2\wedge(s\sigma^1+t\sigma^3)\wedge\sigma^4$ for suitable 
$s,t\in\RR$ with $s^2+t^2\not=0$. Eventually, $[\alpha_7,\gamma']$ is on the 
same $G$-orbit as  
$[\alpha_7,\gamma_0]$.
\qed 

Combining the description of the moduli space given in Proposition \ref{P3} with 
Propositions \ref{P4} -- \ref{Pl} and the computations of $\cHQ_0$ for $\dim 
\fl\le3$ in \cite{KO1} and \cite{KO2} we obtain the following result. We use the 
2-forms $\alpha_1,\dots,\alpha_7,\alpha_5',\alpha_6'$ and the 3-form $\gamma_0$ 
introduced before Proposition \ref{Lh1}. 
\begin{theo} \label{T1}
If $(\fg,\ip)$ is an indecomposable non-abelian nilpotent metric Lie algebra of 
dimension at most 10, then it is 
isomorphic to 
$\fd_{\alpha,\gamma}(\fl,\fa)$ for exactly one of the data in the following 
list: 
\begin{enumerate} 
\item $\fl=\RR^5$ 
\begin{itemize} 
    \item[] $\fa=0$, $\alpha=0$, $\gamma=(\sigma^{1}\wedge\sigma^{2} 
 +\sigma^{3}\wedge\sigma^{4})\wedge\sigma^{5}$; 
\end{itemize}   
\item $\fl=\fg_{5,2} =\{ [X_1,X_2]=Y,\ [X_1,X_3]=Z\}$ 
\begin{enumerate}
\item[] $\fa=0$, $\alpha=0$,  
 $\gamma\in\{\sigma^{1}\wedge\sigma^{Y}\wedge\sigma^{Z},\, 
\sigma^{1}\wedge\sigma^{Y}\wedge\sigma^{Z}+\sigma^{2} 
\wedge\sigma^{3}\wedge\sigma^{Z}\} $; 
\end{enumerate}   
\item $\fl=\fg_{4,1} =\{ [X_1,Z]=Y,\ [X_1,X_2]=Z\}$ 
\begin{enumerate} 
\item $\fa\in \{\RR^{1},\ \RR^{1,0}\}$, $A\in\fa$ fixed unit vector,\\ 
 $\alpha=\sigma^{1}\wedge\sigma^{Y}\otimes A$, \\ 
 $\gamma\in\{0,\,\sigma^{2}\wedge\sigma^{Y}\wedge\sigma^{Z}\}$;
\item $\fa\in \{\RR^{2},\ \RR^{2,0}\}$ with fixed 
 orthonormal basis $A_{1},A_{2}$,\\ 
 $\alpha=\sigma^{1}\wedge\sigma^{Y}\otimes A_{1}+ 
 \sigma^{2}\wedge\sigma^{Z}\otimes A_{2} $,  
 $\gamma=0$;
\item $\fa=\RR^{1,1}$ with fixed 
 orthonormal basis $A_{1},A_{2}$,\\
 $\alpha\in\{\sigma^{1}\wedge\sigma^{Y}\otimes A_{1}+ 
 \sigma^{2}\wedge\sigma^{Z}\otimes A_{2},\ \sigma^{2}\wedge\sigma^{Z}\otimes 
A_{1}+ 
 \sigma^{1}\wedge\sigma^{Y}\otimes A_{2} \}$, \\ 
 $\gamma=0$;  
 \end{enumerate}  
\item $\fl=\fh(1)\oplus\RR^1$ 
\begin{enumerate} 
\item $\fa\in\{\RR^2,\ \RR^{2,0}\}$ with fixed 
 orthonormal basis $A_{1},A_{2}$ \\ 
 $(\alpha,\gamma)\in\{(\alpha_1,0),\,(\alpha_5,0),\,(\alpha_5,\gamma 
_0),\,(\alpha_6,0),\,(\alpha_6,\gamma 
_0)\}$; 
\item$\fa=\RR^{1,1}$, \\ 
$(\alpha,\gamma)\in\{(\alpha_1,0),\,(\alpha_2,0),\,(\alpha_3,0)\}$, 
where $A_{1},A_{2}$ is a fixed Witt basis, or\\ 
$(\alpha,\gamma)\in\{(\alpha_5,0),\,(\alpha_5,\gamma_0),\, 
(\alpha_6,0),\,(\alpha_6,\gamma_0),\,(\alpha_5',0),\,(\alpha_5',\gamma_0),\, 
(\alpha_6',0),\,(\alpha_6',\gamma_0)\}$, where $A_2,A_1$ is a fixed 
orthonormal 
basis 
of $\fa$; 
\item $\fa\in \{\RR^1,\ \RR^{1,0}\}$, \\ 
$(\alpha,\gamma)=(\alpha_7,\gamma_0)$, where $A_{1}$ is a fixed unit vector in 
$\fa$; 
\end{enumerate} 
\item $\fl=\RR^4$ 
\begin{enumerate} 
\item $\fa\in\{\RR^2,\,\RR^{2,0}\}$ with fixed 
 orthonormal basis $A_{1},A_{2}$,\\ 
 $(\alpha,\gamma)\in\{(\alpha_1,0),\,(\alpha_4,\sigma^1\wedge\sigma^2        
\wedge\sigma^4)\}$; 
\item$\fa=\RR^{1,1}$ with fixed 
 Witt basis $A_{1},A_{2}$,\\ 
$(\alpha,\gamma)\in\{(\alpha_1,0), \,(\alpha_2,0),\,(\alpha_3,0),\, 
(\alpha_4,\sigma^1\wedge\sigma^2\wedge\sigma^4)\}$; 
\item $\fa\in\{\RR^1,\ \RR^{1,0}\}$,\\ 
$(\alpha,\gamma)=(\alpha_7,\gamma_0)$, 
where $A_1$ is a fixed unit vector in 
$\fa$; 
\end{enumerate} 
\item $\fl=\fh(1)=\{[X_1,X_2]=Y\}$ 
\begin{enumerate} 
\item $\fa\in \{\RR^{1},\ \RR^{1,0}\}$,\\ 
 $\alpha=\sigma^{1}\wedge\sigma^{Y}\otimes A$, where $A$ is a 
 fixed unit vector in $\fa$,\\ 
 $\gamma=0$; 
\item $\fa\in \{\RR^{2},\ \RR^{2,0},\ \RR^{1,1}\}$ with fixed 
 orthonormal basis $A_{1},A_{2}$,\\ 
 $\alpha=\sigma^{1}\wedge\sigma^{Y}\otimes A_{1}+ 
 \sigma^{2}\wedge\sigma^{Y}\otimes A_{2}$, \\ 
 $\gamma=0$; 
\end{enumerate} 
\item $\fl=\RR^3$ 
\begin{enumerate} 
\item $\fa=0$, $\alpha=0$, 
 $\gamma=\sigma^{1}\wedge\sigma^{2}\wedge\sigma^{3}$;  
\item $\fa\in \{\RR^{2},\ \RR^{2,0},\ \RR^{1,1}\}$ with fixed 
 orthonormal basis $A_{1},A_{2}$,\\ 
 $\alpha=\sigma^{1}\wedge\sigma^{2}\otimes A_{1}+ 
 \sigma^{1}\wedge\sigma^{3}\otimes A_{2}$,\\ 
 $\gamma=0$; 
\item $\fa\in \{\RR^{3},\ \RR^{2,1},\ \RR^{1,2},\ \RR^{3,0}\}$ with fixed 
 orthonormal basis $A_{1},A_{2},A_{3}$,\\ 
 $\alpha=\sigma^{1}\wedge\sigma^{2}\otimes A_{1}+ 
 \sigma^{1}\wedge\sigma^{3}\otimes A_{2}+ 
 \sigma^{2}\wedge\sigma^{3}\otimes A_{3}$, \\ 
 $\gamma=0$; 
\end{enumerate} 
\item $\fl=\RR^2$ 
\begin{itemize} 
\item[] $\fa\in \{\RR^{1},\ \RR^{1,0}\}$,\\ 
 $\alpha=\sigma^{1}\wedge\sigma^{2}\otimes A$, where $A$ is a 
 fixed unit vector in $\fa$,\\ 
 $\gamma=0$. 
\end{itemize} 
\end{enumerate} 
\end{theo} 

{\bf Acknowledgement} I would like to thank my friend Martin Olbrich for all his 
valuable comments and for the mathematical corrections which he brought to my 
attention. 

Moreover, I would like to thank Yves Cornulier for sending me a preprint of his paper {\sl On the Koszul map of Lie algebras}. Motivated by this paper, I checked the list of metric nilpotent Lie algebras of dimension $\le 10$ that I gave in a former version  and realised that it contains a mistake in case 3. Independently, Yves Cornulier and Louis Magnin found isomorphisms that exist between some of the metric Lie algebras in case 3 of the former list.

\vspace{0.8cm} 
{\footnotesize 

Ines Kath\\ 
Institut f\"ur Mathematik und Informatik\\ 
Ernst-Moritz-Arndt-Universit\"at Greifswald\\ 
Walther-Rathenau-Str. 47\\
email: ines.kath@uni-greifswald.de} 

\end{document}